# SOME SHARP PERFORMANCE BOUNDS FOR LEAST SQUARES REGRESSION WITH $L_1$ REGULARIZATION

By Tong Zhang[1]

*Rutgers University*

We derive sharp performance bounds for least squares regression with $L_1$ regularization from parameter estimation accuracy and feature selection quality perspectives. The main result proved for $L_1$ regularization extends a similar result in [*Ann. Statist.* **35** (2007) 2313–2351] for the Dantzig selector. It gives an affirmative answer to an open question in [*Ann. Statist.* **35** (2007) 2358–2364]. Moreover, the result leads to an extended view of feature selection that allows less restrictive conditions than some recent work. Based on the theoretical insights, a novel two-stage $L_1$-regularization procedure with selective penalization is analyzed. It is shown that if the target parameter vector can be decomposed as the sum of a sparse parameter vector with large coefficients and another less sparse vector with relatively small coefficients, then the two-stage procedure can lead to improved performance.

**1. Introduction.** Consider a set of input vectors $\mathbf{x}_1, \ldots, \mathbf{x}_n \in R^d$ with corresponding desired output variables $\mathbf{y}_1, \ldots, \mathbf{y}_n$. We use $d$ instead of the more conventional $p$ to denote data dimensionality, because the symbol $p$ is used for another purpose. The task of supervised learning is to estimate the functional relationship $y \approx f(\mathbf{x})$ between the input $\mathbf{x}$ and the output variable $y$ from the training examples $\{(\mathbf{x}_1, \mathbf{y}_1), \ldots, (\mathbf{x}_n, \mathbf{y}_n)\}$.

In this paper, we consider the linear prediction model $f(\mathbf{x}) = \beta^T \mathbf{x}$ and focus on least squares for simplicity. A commonly used estimation method is $L_1$-regularized empirical risk minimization (aka Lasso)

$$(1) \qquad \hat{\beta} = \arg\min_{\beta \in R^d} \left[ \frac{1}{n} \sum_{i=1}^{n} (\beta^T \mathbf{x}_i - \mathbf{y}_i)^2 + \lambda \|\beta\|_1 \right],$$

Received June 2008; revised September 2008.

[1]Supported in part by NSF Grant DMS-07-06805.

AMS 2000 subject classifications. Primary 62G05; secondary 62J05.

*Key words and phrases.* $L_1$ regularization, Lasso, regression, sparsity, variable selection, parameter estimation.

where $\lambda \geq 0$ is an appropriate regularization parameter.

We are specifically interested in two related themes: parameter estimation accuracy and feature selection quality. A general convergence theorem is established in Section 4, which has two consequences. First, the theorem implies a parameter estimation accuracy result for standard Lasso that extends the main result of Dantzig selector in [4]. A detailed comparison is given in Section 6. This result provides an affirmative answer to an open question in [10] concerning whether a bound similar to that of [4] holds for Lasso. Second, we show, in Section 7, that the general theorem in Section 4 can be used to study the feature selection quality of Lasso. In this context, we consider an extended view of feature selection by selecting features with estimated coefficients larger than a certain nonzero threshold. This method is different from [20], which only considered zero threshold. An interesting consequence of our method is the consistency of feature selection, even when the irrepresentable condition of [20] (the condition is necessary in their approach) is violated. Moreover, the combination of our parameter estimation and feature selection results suggest that the standard Lasso might be sub-optimal when the target can be decomposed as a sparse parameter vector with large coefficients, plus another less sparse vector with small coefficients. A two-stage selective penalization procedure is proposed in Section 8 to remedy the problem. We obtain a parameter estimation accuracy result for this procedure that can improve the corresponding result of the standard (one-stage) Lasso under appropriate conditions.

For simplicity, most results (except for Theorem 4.1) in this paper are stated under the fixed design situation (i.e., $\mathbf{x}_i$ are fixed, while $\mathbf{y}_i$ are random). However, with small modifications, they can also be applied to random design.

## 2. Related work.
In the literature, there are typically three types of results for learning a sparse approximate target vector $\bar{\beta} = [\bar{\beta}_1, \dots, \bar{\beta}_d] \in R^d$ such that $\mathbf{E}(\mathbf{y}|\mathbf{x}) \approx \bar{\beta}^T \mathbf{x}$. These results are as follows:

1. Feature selection accuracy. Identify nonzero coefficients (e.g., [5, 18, 19, 20]), or more generally, identify features with target coefficients larger than a certain threshold (see Section 7). That is, we are interested in identifying the relevant feature set $\{j : |\bar{\beta}_j| > \alpha\}$ for some threshold $\alpha \geq 0$.
2. Parameter estimation accuracy. How accurate is the estimated parameter, comparing to the approximate target $\bar{\beta}$, measured in a certain norm (e.g., [1, 2, 4, 5, 9, 13, 17, 19])? That is, let $\hat{\beta}$ be the estimated parameter; we are interested in developing a bound for $\|\hat{\beta} - \bar{\beta}\|_p$ for some $p$. Theorems 4.1 and 8.1 give such results.
3. Prediction accuracy. The prediction performance of the estimated parameter, both in fixed and random design settings (e.g., [2, 3, 5, 11, 13, 14]). For example, in fixed design, we are interested in a bound for $\frac{1}{n}\sum_{i=1}^{n}(\hat{\beta}^T\mathbf{x}_i - \mathbf{E}\mathbf{y}_i)^2$ or the related quantity $\frac{1}{n}\sum_{i=1}^{n}((\hat{\beta}-\bar{\beta})^T\mathbf{x}_i)^2$.

In general, good feature selection implies good parameter estimation, and good parameter estimation implies good prediction accuracy. However, the reverse directions do not usually hold. In this paper, we focus on the first two aspects, feature selection and parameter estimation, as well as their inter-relationship. Due to the space limitation, the prediction accuracy of Lasso is not consider in this paper. However, it is a relatively straight-forward consequence of parameter estimation bounds with $p=1$ and $p=2$.

As mentioned in the Introduction, one motivation of this paper is to develop a parameter estimation bound for $L_1$ regularization directly comparable to that of the Dantzig selector in [4]. Compared to [4], where a parameter estimation error bound in 2-norm is proved for the Dantzig selector, Theorem 4.1 in Section 4 presents a more general bound for $p$-norm where $p \in [1, \infty]$. We are particularly interested in a bound in $\infty$-norm because such a bound immediately induces a result on the feature selection accuracy of Lasso. This point of view is taken in Section 7, where feature selection is considered. Achieving good feature selection is important, because it can be used to improve the standard one-stage Lasso. In Section 8, we develop this observation and show that a two-stage method with good feature selection achieves a bound better than that of one-stage Lasso. Experiments in Section 9 confirm this theoretical observation.

Since the development of this paper relies heavily on parameter estimation accuracy of Lasso in $p$-norm, it is different and complements earlier work on Lasso given above. Among earlier work, prediction accuracy bounds for Lasso were derived in [2, 3] and [5] under mutual incoherence conditions (introduced in [9]) that are generally regarded as stronger than the sparse eigenvalue conditions employed by [4]. This is because it is easier for a random matrix to satisfy sparse eigenvalue conditions than mutual incoherence conditions. A more detailed discussion on this point is given at the end of Section 4. Moreover, the relationship of different quantities are presented in Section 3. As we shall see from Section 3, mutual incoherence conditions are also stronger than conditions required for deriving $p$-norm estimation error bounds in this paper. Therefore, under appropriate mutual incoherence conditions, our analysis leads to sharp $p$-norm parameter estimation bounds for all $p \in [1, \infty]$ in Corollary 4.1. In comparison, sharp $p$-norm parameter estimation bounds cannot be directed derived from prediction error bounds studied in some earlier work.

We shall also point out that some 1-norm parameter estimation bounds were established in [2], but not for $p > 1$. At the same time this paper was written, related parameter estimation error bounds were also obtained in [1], under appropriate sparse eigenvalue conditions, both for the Dantzig selector and for Lasso. However, the results are only for $p \in [0, 2]$ and only for truly sparse targets [i.e., $\mathbf{E}(\mathbf{y}|\mathbf{x}) = \bar{\beta}^T \mathbf{x}$ for some sparse $\bar{\beta}$]. In particular, their parameter estimation, bound with $p = 2$, does not reproduce the main result of [4], while our bounds in this paper (which allows approximate sparse target) do. We shall point out that that, in addition to parameter estimation error bounds, a prediction error bound that allows approximately sparse target was also obtained, in [1], in a form quite similar to the parameter estimation bounds of Lasso in this paper and that of Dantzig selector in [4]. However, that result does not immediately imply a bound on $p$-norm parameter estimation error. A similar bound was derived in [5] for Lasso but not as elaborated as that in [4]. Another related work is [19], which contains a 2-norm parameter estimation error bound for Lasso but has a cruder form than ours. In particular, their result is worse than the form given in [4], as well as the first claim of Theorem 4.1 in this paper. A similar 2-norm parameter estimation error bound, but only for truly sparse targets, can be found in [17]. In [9], the authors derived 2-norm estimation bound for Lasso with approximate sparse targets under mutual incoherence conditions but without stochastic noise. We shall note that the "noise" in their paper is not random and corresponds to approximation error in our notation, as discussed in Section 5. Their result is thus weaker than our result in Corollary 5.1. In addition to the above work, prediction error bounds were also obtained in [13, 14] and [11] for general loss functions and random design. However, such results cannot be used to derive $p$-norm parameter estimation bound, which we consider here.

**3. Conditions on design matrix.** In order to obtain good bounds, it is necessary to impose conditions on the design matrix that generally specifies that small diagonal blocks of the design matrix are nonsingular. For example, mutual incoherence conditions [9] or sparse eigenvalue conditions [19]. The sparse eigenvalue condition is also known as RIP (restricted isometry property) in the compressed sensing literature, which was first introduced in [7].

Mutual incoherence conditions are usually more restrictive than sparse eigenvalue conditions. That is, a matrix that satisfies an appropriate mutual incoherence condition will also satisfy the necessary sparse eigenvalue condition in our analysis, but the reverse direction does not hold. For example, we will see from the discussion at the end of Section 4 that, for random design matrices, more samples are needed in order to satisfy the mutual incoherence condition than to satisfy the sparse eigenvalue condition. In our analysis, the weaker sparse eigenvalue condition can be used to obtain sharp bounds for 2-norm parameter estimation error. Since we are interested in general $p$-norm parameter estimation error bounds, other conditions on the design matrix will be considered. They can be regarded as generalizations of the sparse eigenvalue condition, or RIP. All of these conditions are weaker than the mutual incoherence condition.

We introduce the following definitions that specify properties of submatrices of a large matrix $A$. These quantities (when used with the design matrix $\hat{A} = \frac{1}{n} \sum_{i=1}^{n} \mathbf{x}_i \mathbf{x}_i^T$) will appear in our result.

DEFINITION 3.1. The $p$-norm of a vector $\beta = [\beta_1, \ldots, \beta_d] \in R^d$ is defined as $\|\beta\|_p = (\sum_{j=1}^{d} |\beta_j|^p)^{1/p}$. Given a positive semi-definite matrix $A \in R^{d \times d}$, and given $\ell, k \geq 1$ such that $\ell + k \leq d$, let $I, J$ be disjoint subsets of $\{1, \ldots, d\}$ with $k$ and $\ell$ elements, respectively. Also, let $A_{I,I} \in R^{k \times k}$ be the restriction of $A$ to indices $I$, $A_{I,J} \in R^{k \times \ell}$ be the restriction of $A$ to indices $I$ on the left and $J$ on the right. Define, for $p \in [1, \infty]$,

$$\rho_{A,k}^{(p)} = \sup_{\mathbf{v} \in R^k, I} \frac{\|A_{I,I}\mathbf{v}\|_p}{\|\mathbf{v}\|_p}, \qquad \theta_{A,k,\ell}^{(p)} = \sup_{\mathbf{u} \in R^\ell, I, J} \frac{\|A_{I,J}\mathbf{u}\|_p}{\|\mathbf{u}\|_\infty},$$

$$\mu_{A,k}^{(p)} = \inf_{\mathbf{v} \in R^k, I} \frac{\|A_{I,I}\mathbf{v}\|_p}{\|\mathbf{v}\|_p}, \qquad \gamma_{A,k,\ell}^{(p)} = \sup_{\mathbf{u} \in R^\ell, I, J} \frac{\|A_{I,I}^{-1} A_{I,J}\mathbf{u}\|_p}{\|\mathbf{u}\|_\infty}.$$

Moreover, for all $\mathbf{v} = [\mathbf{v}_1, \ldots, \mathbf{v}_k] \in R^k$, define $\mathbf{v}^{p-1} = [|\mathbf{v}_1|^{p-1} \operatorname{sgn}(\mathbf{v}_1), \ldots, |\mathbf{v}_k|^{p-1} \operatorname{sgn}(\mathbf{v}_k)]$, and

$$\omega_{A,k}^{(p)} = \inf_{\mathbf{v} \in R^k, I} \frac{\max(0, \mathbf{v}^T A_{I,I} \mathbf{v}^{p-1})}{\|\mathbf{v}\|_p^p},$$

$$\pi_{A,k,\ell}^{(p)} = \sup_{\mathbf{v} \in R^k, \mathbf{u} \in R^\ell, I, J} \frac{(\mathbf{v}^{p-1})^T A_{I,J}\mathbf{u} \|\mathbf{v}\|_p}{\max(0, \mathbf{v}^T A_{I,I} \mathbf{v}^{p-1}) \|\mathbf{u}\|_\infty}.$$

The ratio $\rho_{A,k}^{(p)} / \mu_{A,k}^{(p)}$ measures the closeness to the identity matrix of $k \times k$ diagonal sub-matrices of $A$. The RIP concept in [7] can be regarded as $\rho_{A,k}^{(2)} / \mu_{A,k}^{(2)}$ in our notation. The quantities $\theta_{A,k,\ell}^{(p)}$, $\gamma_{A,k,\ell}^{(p)}$ and $\pi_{A,k,\ell}^{(p)}$ measures the closeness to zero of the $k \times \ell$ off diagonal blocks of $A$. Note that $\mu_{A,k}^{(2)} = \omega_{A,k}^{(2)}$ and $\rho_{A,k}^{(2)}$ are the smallest and largest eigenvalues of $k \times k$ diagonal blocks of $A$. It is easy to see that the inequalities $\mu_{A,k}^{(p)} \leq \rho_{A,k}^{(p)}$ hold. Moreover, we can also obtain bounds on $\theta_{A,k,\ell}^{(p)}$, $\gamma_{A,k,\ell}^{(p)}$ and $\pi_{A,k,\ell}^{(2)}$ using eigenvalues of sub-matrices of $A$. The bounds essentially say that if diagonal sub-matrices of $A$ of size $k + \ell$ are well-conditioned, then the quantities $\theta_{A,k,\ell}^{(2)}$, $\gamma_{A,k,\ell}^{(2)}$ and $\pi_{A,k,\ell}^{(2)}$ are $O(\sqrt{\ell})$.

PROPOSITION 3.1. *The following inequalities hold:*

$$\theta_{A,k,\ell}^{(2)} \leq \ell^{1/2} \sqrt{(\rho_{A,k}^{(2)} - \mu_{A,\ell+k}^{(2)})(\rho_{A,\ell}^{(2)} - \mu_{A,\ell+k}^{(2)})},$$

$$\theta_{A,k,\ell}^{(p)} \leq k^{\max(0,1/p-0.5)} \theta_{A,k,\ell}^{(2)},$$

$$\pi_{A,k,\ell}^{(2)} \leq \frac{\ell^{1/2}}{2} \sqrt{\rho_{A,\ell}^{(2)}/\mu_{A,k+\ell}^{(2)} - 1}, \qquad \pi_{A,k,\ell}^{(p)} \leq \theta_{A,k,\ell}^{(p)}/\omega_{A,k}^{(p)},$$

$$\gamma_{A,k,\ell}^{(p)} \leq k^{\max(0,1/p-0.5)} \gamma_{A,k,\ell}^{(2)}, \qquad \gamma_{A,k,\ell}^{(p)} \leq \theta_{A,k,\ell}^{(p)}/\mu_{A,k}^{(p)},$$

$$\min_i A_{i,i} - \sup_I \|A_{I,I} - \operatorname{diag}(A_{I,I})\|_p \leq \omega_{A,k}^{(p)} \leq \mu_{A,k}^{(p)},$$

*where, for a matrix $B$, $\operatorname{diag}(B)$ is the diagonal of $B$, and $\|B\|_p = \sup_{\mathbf{u}}(\|B\mathbf{u}\|_p/\|\mathbf{u}\|_p)$.*

The last inequality in the above proposition shows that $\mu_{A,k}^{(p)} > 0$ and $\omega_{A,k}^{(p)} > 0$ when $A$ has a certain diagonal dominance (in $p$-norm) property for its $k \times k$ blocks.

Finally, we state a result that bounds all quantities defined here using the mutual incoherence concept of [9]. This result shows that mutual incoherence is a stronger notation than all quantities we employ in this paper. Although more complicated, by using these less restrictive quantities in Theorem 4.1, we obtain stronger results than using the mutual incoherence condition (Corollary 4.1). For simplicity, we consider diagonally normalized $A$ such that $A_{i,i} = 1$ for all $i$.

PROPOSITION 3.2.    *Given a matrix $A \in R^{d \times d}$ and assuming that $A_{i,i} = 1$ for all $i$, define the mutual coherence coefficient as $M_A = \sup_{i \neq j} |A_{i,j}|$. Then the following bounds hold:*

- $\rho_{A,k}^{(p)} \leq 1 + M_A k;$
- $\mu_{A,k}^{(p)} \geq \omega_{A,k}^{(p)} \geq 1 - M_A k;$
- $\theta_{A,k,\ell}^{(p)} \leq M_A k^{1/p} \ell;$
- $\pi_{A,k,\ell}^{(p)} \leq M_A k^{1/p} \ell / \max(0, 1 - M_A k);$
- $\gamma_{A,k,\ell}^{(p)} \leq M_A k^{1/p} \ell / \max(0, 1 - M_A k).$

## 4. A general performance bound for $L_1$ regularization.

For simplicity, we assume sub-Gaussian noise as follows. We use $\mathbf{x}_{i,j}$ to indicate the $j$th component of vector $\mathbf{x}_i \in R^d$.

ASSUMPTION 4.1.    *Assume that, conditioned on $\{\mathbf{x}_i\}_{i=1,\ldots,n}$, $\{\mathbf{y}_i\}_{i=1,\ldots,n}$ are independent (but not necessarily identically distributed) sub-Gaussians. There exists a constant $\sigma \geq 0$ such that $\forall i$ and $\forall t \in R$,*

$$\mathbf{E}_{\mathbf{y}_i} e^{t(\mathbf{y}_i - \mathbf{E}\mathbf{y}_i)} | \{\mathbf{x}_i\}_{i=1,\ldots,n} \leq e^{\sigma^2 t^2/2}.$$

Both Gaussian and bounded random variables are sub-Gaussian using the above definition. For example, if a random variable $\xi \in [a, b]$, then $\mathbf{E}_\xi e^{t(\xi - \mathbf{E}\xi)} \leq e^{(b-a)^2 t^2 / 8}$. If a random variable is Gaussian, $\xi \sim N(0, \sigma^2)$, then $\mathbf{E}_\xi e^{t\xi} \leq e^{\sigma^2 t^2 / 2}$.

For convenience, we also introduce the following definition.

DEFINITION 4.1. Let $\beta = [\beta_1, \ldots, \beta_d] \in R^d$ and $\alpha \geq 0$. We define the set of relevant features with threshold $\alpha$ as

$$\mathrm{supp}_\alpha(\beta) = \{j : |\beta_j| > \alpha\}.$$

Moreover, if $|\beta_{(1)}| \geq \cdots \geq |\beta_{(d)}|$ are in descending order, then define

$$r_k^{(p)}(\beta) = \left( \sum_{j=k+1}^d |\beta_{(j)}|^p \right)^{1/p}$$

as the $p$-norm of the $d - k$ smallest components (in absolute value) of $\beta$.

Consider a target parameter vector $\bar{\beta} \in R^d$ that is approximately sparse. Note that we do not assume that $\bar{\beta}$ is the true target; that is, $\bar{\beta}^T \mathbf{x}_i$ may not equal to $\mathbf{E} y_i$. We only assume that this holds approximately, and we are interested in how well we can estimate $\bar{\beta}$ using (1). In particular, the approximation quality (or its closeness to the true model) of any $\bar{\beta}$ is measured by $\|\frac{1}{n} \sum_{i=1}^n (\bar{\beta}^T \mathbf{x}_i - \mathbf{E} y_i) \mathbf{x}_i\|_\infty$ in our analysis. If $\bar{\beta}^T \mathbf{x}_i = \mathbf{E} y_i$, then the underlying model is $y_i = \bar{\beta}^T \mathbf{x}_i + \varepsilon_i$, where $\varepsilon_i$ are independent sub-Gaussian noises. In the more general case, we only need to assume that $\|\frac{1}{n} \sum_{i=1}^n (\bar{\beta}^T \mathbf{x}_i - \mathbf{E} y_i) \mathbf{x}_i\|_\infty$ is small for some approximate target $\bar{\beta}$. The relationship of this quantity and least squares approximation error is discussed in Section 5.

The following theorem holds both in fixed design and in random design. The only difference is that, in the fixed design situation, we may let $a = (\sup_j \hat{A}_{j,j})^{1/2}$, and the condition $(\sup_j \hat{A}_{j,j})^{1/2} \leq a$ automatically holds. In order to simplify the claims, our later results are all stated under the fixed design assumption. In the following theorem, the statement of "with probability $1 - \delta$: if $X$ then $Y$" can also be interpreted as "with probability $1 - \delta$: either $X$ is false or $Y$ is true." We note that, in practice, the condition of the theorem can be combinatorially hard to check, since computing the quantities in Definition 3.1 requires searching over sets of fixed cardinality.

THEOREM 4.1. Let Assumption 4.1 hold, and let $\hat{A} = \frac{1}{n} \sum_{i=1}^n \mathbf{x}_i \mathbf{x}_i^T$. Let $\hat{\beta}$ be the solution of (1). Consider any fixed target vector $\bar{\beta} \in R^d$ and a positive constant $a > 0$. Given $\delta \in (0, 1)$, then, with probability larger than $1 - \delta$, the following two claims hold for $q = 1, p$, and all $k, \ell$ such that $k \leq \ell \leq (d - k)/2$, $t \in (0, 1)$, $p \in [1, \infty]$:

- If $t \leq 1 - \pi_{\hat{A},k+\ell,\ell}^{(p)} k^{1-1/p}/\ell$, $\lambda \geq \frac{4(2-t)}{t}(\sigma a \sqrt{\frac{2}{n}\ln(2d/\delta)} + \|\frac{1}{n}\sum_{i=1}^{n}(\bar{\beta}^T\mathbf{x}_i - \mathbf{E}\mathbf{y}_i)\mathbf{x}_i\|_{\infty})$, and $(\sup_j \hat{A}_{j,j})^{1/2} \leq a$, then

$$\|\hat{\beta} - \bar{\beta}\|_q \leq \frac{8k^{1/q-1/p}}{t\omega_{\hat{A},k+\ell}^{(p)}}[\rho_{\hat{A},k+\ell}^{(p)} r_k^{(p)}(\bar{\beta}) + k^{1/p}\lambda]$$

$$+ \frac{32k^{1/q-1/p}}{t}\pi_{\hat{A},k+\ell,\ell}^{(p)} r_k^{(1)}(\bar{\beta})\ell^{-1}$$

$$+ 4k^{1/q-1/p} r_k^{(p)}(\bar{\beta}) + 4r_k^{(1)}(\bar{\beta})\ell^{1/q-1};$$

- If $t \leq 1 - \gamma_{\hat{A},k+\ell,\ell}^{(p)} k^{1-1/p}/\ell$, $\lambda \geq \frac{4(2-t)}{t}(\sigma a \sqrt{\frac{2}{n}\ln(2d/\delta)} + \|\frac{1}{n}\sum_{i=1}^{n}(\bar{\beta}^T\mathbf{x}_i - \mathbf{E}\mathbf{y}_i)\mathbf{x}_i\|_{\infty})$, and $(\sup_j \hat{A}_{j,j})^{1/2} \leq a$, then

$$\|\hat{\beta} - \bar{\beta}\|_q \leq \frac{8k^{1/q-1/p}}{t}[4\gamma_{\hat{A},k+\ell,\ell}^{(p)}\ell^{-1} r_k^{(1)}(\bar{\beta}) + \lambda(k+\ell)^{1/p}/\mu_{\hat{A},k+\ell}^{(p)}]$$

$$+ 4r_k^{(1)}(\bar{\beta})\ell^{1/q-1}.$$

Although, for a fixed $p$, the above theorem only gives bounds for $\|\hat{\beta} - \bar{\beta}\|_q$ with $q = 1, p$, this information is sufficient to obtain bounds for general $q \in [1,\infty]$. If $q \in [1,p]$, we can use the following interpolation rule (which follows from the Hölder's inequality):

$$\|\Delta\hat{\beta}\|_q \leq \|\Delta\hat{\beta}\|_1^{(p-q)/(qp-q)} \|\Delta\hat{\beta}\|_p^{(pq-p)/(pq-q)}$$

and if $q > p$, we use $\|\Delta\hat{\beta}\|_q \leq \|\Delta\hat{\beta}\|_p$. Although the estimate we obtain when $q \neq p$ is typically worse than the bound achieved at $p = q$ (assuming that the condition of the theorem can be satisfied at $p = q$), it may still be useful, because the condition for the theorem to apply may be easier to satisfy for certain $p$.

It is important to note that the first claim is more refined than the second claim, as it replaces the explicit $\ell$-dependent term $O(\ell^{1/p}\lambda)$ by the term $O(r_k^{(p)}(\bar{\beta}))$, which does not explicitly depend on $\ell$. In order to optimize the bound, we can choose $k = |\operatorname{supp}_{\lambda}(\bar{\beta})|$, which implies that $O(r_k^{(p)}(\bar{\beta})) = O(\ell^{1/p}\lambda + r_{k+\ell}^{(p)}(\bar{\beta})) = O(\ell^{1/p}\lambda + \ell^{1/p-1}r_k^{(1)}(\bar{\beta}))$. This quantity is dominated by the bound in the second claim. However, if $O(r_k^{(p)}(\bar{\beta}))$ is small, then the first claim is much better when $\ell$ is large.

The theorem as stated is not intuitive. In order to obtain a more intuitive bound from the first claim of the theorem, we consider a special case with mutual incoherence condition. The following corollary is a simple consequence of the first claim of Theorem 4.1 (with $q = p$) and Proposition 3.2.

The result shows that mutual incoherence condition is a stronger assumption than the quantities that appear in our analysis.

COROLLARY 4.1. *Let Assumption 4.1 hold, and let* $\hat{A} = \frac{1}{n}\sum_{i=1}^{n}\mathbf{x}_i\mathbf{x}_i^T$, *and assume that* $\hat{A}_{j,j} = 1$ *for all* $j$. *Define* $M_{\hat{A}} = \sup_{i\neq j}|\hat{A}_{i,j}|$. *Let* $\hat{\beta}$ *be the solution of (1). Consider any fixed target vector* $\bar{\beta} \in R^d$. *Given* $\delta \in (0,1)$, *then, with probability larger than* $1 - \delta$, *the following claim holds for all* $k \leq \ell \leq (d-k)/2$, $t \in (0,1)$, $p \in [1,\infty]$: *if* $M_{\hat{A}}(k+\ell) \leq (1-t)/(2-t)$ *and* $\lambda \geq \frac{4(2-t)}{t}(\sigma\sqrt{\frac{2}{n}\ln(2d/\delta)} + \|\frac{1}{n}\sum_{i=1}^{n}(\bar{\beta}^T\mathbf{x}_i - \mathbf{E}\mathbf{y}_i)\mathbf{x}_i\|_\infty)$, *then*

$$\|\hat{\beta} - \bar{\beta}\|_p \leq \frac{8(2-t)}{t}[1.5r_k^{(p)}(\bar{\beta}) + k^{1/p}\lambda] + 4r_k^{(p)}(\bar{\beta}) + \frac{4(8-7t)}{t}r_k^{(1)}(\bar{\beta})\ell^{1/p-1}.$$

The above result is of the form

$$\|\hat{\beta} - \bar{\beta}\|_p = O(k^{1/p}\lambda + r_k^{(p)}(\bar{\beta}) + r_k^{(1)}(\bar{\beta})\ell^{1/p-1}), \tag{2}$$

where we can let $k = |\text{supp}_\lambda(\bar{\beta})|$, so that $k$ is the number of components such that $|\bar{\beta}_j| > \lambda$. Although mutual incoherence is assumed for simplicity, a similar bound holds for any $p$ if we assume that $\hat{A}$ is $p$-diagonal dominant at block size $k + \ell$. Such an assumption is weaker than mutual-incoherence.

To our knowledge, none of the earlier work on Lasso obtained parameter estimation bounds in the form of (2). The first two terms in the bound are what we shall expect from the $L_1$-regularization method (1) and, thus, unlikely to be significantly improved (except for the constants). The first term is the variance term, and the second term is needed because $L_1$ regularization tends to shrink coefficients $j \notin \text{supp}_\lambda(\bar{\beta})$ to zero. Although it is not clear whether the third term can be improved, we shall note that it becomes small if we can choose a large $\ell$. Note that, if $\bar{\beta}$ is the true parameter: $\bar{\beta}^T\mathbf{x}_i = \mathbf{E}\mathbf{y}_i$, then we may take $\lambda = 4(2-t)t^{-1}\sigma\sqrt{\ln(d/\delta)/n}$ in (2). The bound in [4] has a similar form (but with $p=2$), which we will compare in the next section.

Note that, in Theorem 4.1, one can always take $\lambda$ sufficiently large, so that the condition for $\lambda$ is satisfied. Therefore, in order to apply the theorem, one needs either the condition $0 < t \leq 1 - \pi_{\hat{A},k+\ell,\ell}^{(p)}k^{1-1/p}/\ell$ or the condition $0 < t \leq 1 - \gamma_{\hat{A},k+\ell,\ell}^{(p)}k^{1-1/p}/\ell$. They require that small diagonal blocks of $\hat{A}$ are not nearly singular. As pointed out after Proposition 3.1, the condition for the first claim is typically harder to satisfy. For example, as discussed below, even when $p \neq 2$, the requirement $\gamma_{\hat{A},k+\ell,\ell}^{(p)}k^{1-1/p}/\ell < 1$ can always be satisfied when diagonal sub-blocks of $\hat{A}$ at certain size satisfy some eigenvalue conditions, while this is not true for the condition $\pi_{\hat{A},k+\ell,\ell}^{(p)}k^{1-1/p}/\ell < 1$.

In the case of $p = 2$, the condition $0 < 1 - \pi_{\hat{A},k+\ell,\ell}^{(2)}k^{0.5}/\ell$ can always be satisfied if the small diagonal blocks of $\hat{A}$ have eigenvalues bounded from above and below (away from zero). We shall refer to such a condition as *sparse eigenvalue condition* (also see [1, 19]). Indeed, Proposition 3.1 implies that $\pi^{(2)}_{\hat{A},k+\ell,\ell} k^{0.5}/\ell \le 0.5(k/\ell)^{0.5}\sqrt{\rho^{(2)}_{A,\ell}/\mu^{(2)}_{A,k+2\ell} - 1}$. Therefore, this condition can be satisfied if we can find $\ell$ such that $\rho^{(2)}_{A,\ell}/\mu^{(2)}_{A,k+2\ell} \le \ell/k$. In particular, if $\rho^{(2)}_{A,\ell}/\mu^{(2)}_{A,k+2\ell} \le c$ for a constant $c > 0$ when $\ell \le ck$ (which is what we will mean by sparse eigenvalue condition in later discussions), then one can simply take $\ell = ck$.

For $p > 2$, a similar claim holds for the condition $0 < 1 - \gamma^{(p)}_{\hat{A},k+\ell,\ell} k^{1-1/p}/\ell$. Proposition 3.1 implies that $\gamma^{(p)}_{\hat{A},k+\ell,\ell} \le \gamma^{(2)}_{\hat{A},k+\ell,\ell} \le \theta^{(2)}_{\hat{A},k+\ell,\ell} \le \theta^{(2)}_{\hat{A},k+\ell,\ell}/\mu^{(2)}_{\hat{A},k+\ell} \le \sqrt{\ell}\rho^{(2)}_{\hat{A},k+\ell}/\mu^{(2)}_{\hat{A},k+\ell}$. Sparse eigenvalue condition (bounded eigenvalue) at block size of order $k^{2-2/p}$ implies that the condition $0 < 1 - \gamma^{(p)}_{\hat{A},k+\ell,\ell} k^{1-1/p}/\ell$ can be satisfied with an appropriate choice of $\ell = O(k^{2-2/p})$. Under assumptions that are stronger than the sparse eigenvalue condition, one can obtain better and simpler results. For example, this is demonstrated in Corollary 4.1 under the mutual incoherence condition.

Finally, in order to concretely compare the condition of Theorem 4.1 for different $p$, we consider the random design situation (with random vectors $\mathbf{x}_i$), where each component $\mathbf{x}_{i,j}$ is independently drawn from the standard Gaussian distribution $N(0,1)$ ($i = 1, \ldots, n$ and $j = 1, \ldots, d$). This situation is investigated in the compressed sensing literature, such as [7]. In particular, it was shown that, with large probability, the following RIP inequality holds for some constant $c > 0$ ($s \le n \le d$): $|\rho^{(2)}_{\hat{A},s} - 1| + |\mu^{(2)}_{\hat{A},s} - 1| \le c\sqrt{s \ln d/n}$. Now, for $p \ge 2$, using Proposition 3.1, it is not hard to show that (we shall skip the detailed derivation here because it is not essential to the main point of this paper), for $k \le \ell$, $\theta^{(p)}_{\hat{A},k+\ell,\ell} \le 4c\ell\sqrt{\ln d/n}$ and $\omega^{(p)}_{\hat{A},k+\ell} \ge 1 - 2c\sqrt{\ell^{2-2/p} \ln d/n}$. Therefore, the condition $0.5 \le 1 - \gamma^{(p)}_{\hat{A},k+\ell,\ell} k^{1-1/p}/\ell$ in Theorem 4.1 holds with large probability, as long as $n \ge 256c^2 \ell^{2-2/p} \ln d$. Therefore, in order to apply the theorem with fixed $k \le \ell$ and $d$, the larger $p$ is, the larger sample size $n$ has to be. In comparison, the mutual incoherence condition of Corollary 4.1 is satisfied when $n \ge c'\ell^2 \ln d$ for some constant $c' > 0$.

## 5. Noise and approximation error.

In Theorem 4.1 and Corollary 4.1, we do not assume that $\bar{\beta}$ is the true parameter that generates $\mathbf{Ey}_i$. The bounds depend on the quantity $\|\frac{1}{n}\sum_{i=1}^{n}(\bar{\beta}^T\mathbf{x}_i - \mathbf{Ey}_i)\mathbf{x}_i\|_\infty$ to measure how close $\bar{\beta}$ is different from the true parameter. This quantity may be regarded as an algebraic definition of noise, in that it behaves like stochastic noise.

The following proposition shows that, if the least squares error achieved by $\bar{\beta}$ (which is often called approximation error) is small, then the algebraic noise level (as defined above) is also small. However, the reverse is not true. For example, for the identity design matrix and $\bar{\beta}_*^T \mathbf{x}_i = \mathbf{E}\mathbf{y}_i$, the algebraic noise is $\|\frac{1}{n}\sum_{i=1}^n (\bar{\beta}^T\mathbf{x}_i - \mathbf{E}\mathbf{y}_i)\mathbf{x}_i\|_\infty = \|\bar{\beta} - \beta_*\|_\infty$, but the least squares approximation error is $\|\bar{\beta} - \beta_*\|_2^2$. Therefore, in general, algebraic noise can be small, even when the least squares approximation error is large.

PROPOSITION 5.1. *Let* $\hat{A} = \frac{1}{n}\sum_{i=1}^n \mathbf{x}_i\mathbf{x}_i^T$ *and* $a = (\sup_j \hat{A}_{j,j})^{1/2}$. *Given* $k \geq 0$, *there exists* $\bar{\beta}^{(k)} \in R^d$ *such that*

$$\left\|\frac{1}{n}\sum_{i=1}^n (\bar{\beta}^{(k)T}\mathbf{x}_i - \mathbf{E}\mathbf{y}_i)\mathbf{x}_i\right\|_\infty \leq \frac{a}{\sqrt{k+1}}\left(n^{-1}\sum_{i=1}^n (\bar{\beta}^T\mathbf{x}_i - \mathbf{E}\mathbf{y}_i)^2\right)^{1/2},$$

$\text{supp}_0(\bar{\beta}^{(k)} - \bar{\beta}) \leq k$, *and* $\|\bar{\beta}^{(k)} - \bar{\beta}\|_2 \leq 2(n^{-1}\sum_{i=1}^n (\bar{\beta}^T\mathbf{x}_i - \mathbf{E}\mathbf{y}_i)^2)^{1/2}/\sqrt{\mu_{\hat{A},k}^{(2)}}.$

This proposition can be combined with Theorem 4.1 or Corollary 4.1 to derive bounds in terms of the approximation error $n^{-1}\sum_{i=1}^n (\bar{\beta}^T\mathbf{x}_i - \mathbf{E}\mathbf{y}_i)^2$ instead of the algebraic noise $\|\frac{1}{n}\sum_{i=1}^n (\bar{\beta}^T\mathbf{x}_i - \mathbf{E}\mathbf{y}_i)\mathbf{x}_i\|_\infty$. For example, as a simple consequence of Corollary 4.1, we have the following bound. A similar but less general result [with $\sigma = 0$ and $|\text{supp}_0(\bar{\beta})| = k$] was presented in [9].

COROLLARY 5.1. *Let Assumption 4.1 hold, let* $\hat{A} = \frac{1}{n}\sum_{i=1}^n \mathbf{x}_i\mathbf{x}_i^T$ *and assume that* $\hat{A}_{j,j} = 1$ *for all* $j$. *Define* $M_{\hat{A}} = \sup_{i \neq j} |\hat{A}_{i,j}|$. *Let* $\hat{\beta}$ *be the solution of* (1). *Consider any fixed target vector* $\bar{\beta} \in R^d$. *Given* $\delta \in (0,1)$, *then, with probability larger than* $1 - \delta$, *the following claim holds for all* $2k \leq \ell \leq (d-2k)/2$, $t \in (0,1)$, $p \in [1,\infty]$. *If* $M_{\hat{A}}(2k+\ell) \leq (1-t)/(2-t)$ *and* $\lambda \geq \frac{4(2-t)}{t}(\sigma\sqrt{\frac{2}{n}\ln(2d/\delta)} + \varepsilon/\sqrt{k+1})$, *then*

$$\|\hat{\beta} - \bar{\beta}\|_2 \leq \frac{8(2-t)}{t}[1.5r_k^{(2)}(\bar{\beta}) + (2k)^{1/2}\lambda] + 4r_k^{(2)}(\bar{\beta})$$

$$+ \frac{4(8-7t)}{t}r_k^{(1)}(\bar{\beta})\ell^{-1/2} + 4\varepsilon,$$

*where* $\varepsilon = (n^{-1}\sum_{i=1}^n (\bar{\beta}^T\mathbf{x}_i - \mathbf{E}\mathbf{y}_i)^2)^{1/2}$.

A similar result holds under the sparse eigenvalue condition. We should point out that, in $L_1$ regularization, the behavior of stochastic noise ($\sigma > 0$) is similar to that of the algebraic noise introduced above, but it is very different from the least squares approximation error $\varepsilon$. In particular, the so-called bias of $L_1$ regularization shows up in the stochastic noise term but not in the least squares approximation error term. If we set $\sigma = 0$ but $\varepsilon \neq 0$, our analysis of the two-stage procedure in Section 8 will not improve that of the standard Lasso given in Corollary 5.1, simply because the two-stage procedure does not improve the term involving the approximation error $\varepsilon$. However, the benefit of the two-stage procedure clearly shows up in the stochastic noise term. For this reason, it is important to distinguish the true stochastic noise and the approximation error $\varepsilon$, and to develop analysis that includes both stochastic noise and approximation error.

**6. Dantzig selector versus $L_1$ regularization.** Recently, Candes and Tao proposed an estimator, called the *Dantzig selector* in [4], and proved a very strong performance bound for this new method. However, it was observed [8, 10] and [12] that the performance of Lasso is comparable to that of the Dantzig selector. Consequently, the authors of [10] asked whether a performance bound similar to the Dantzig selector holds for Lasso as well. In this context, we observe that a simple but important consequence of the first claim of Theorem 4.1 leads to a bound for $L_1$ regularization that reproduces the main result for Dantzig selector in [4]. We restate the result below, which provides an affirmative answer to the above mentioned open question of [10].

COROLLARY 6.1. *Let Assumption 4.1 hold, and let* $\hat{A} = \frac{1}{n}\sum_{i=1}^n \mathbf{x}_i \mathbf{x}_i^T$ *and* $a = (\sup_j \hat{A}_{j,j})^{1/2}$. *Consider the true target vector* $\bar{\beta}$ *such that* $\mathbf{E}y = \bar{\beta}^T \mathbf{x}$. *Define* $\hat{A} = \frac{1}{n}\sum_{i=1}^n \mathbf{x}_i \mathbf{x}_i^T$. *Let* $\hat{\beta}$ *be the solution of (1). Given* $\delta \in (0,1)$, *then, with probability larger than* $1-\delta$, *the following claim holds for all* $(d-k)/2 \geq \ell \geq k$. *If* $t = 1 - \pi_{\hat{A},k+\ell,\ell}^{(2)} k^{0.5}/\ell > 0$, $\lambda \geq 4(2-t)t^{-1}\sigma a\sqrt{2\ln(2d/\delta)/n}$, *then*

$$\|\hat{\beta} - \bar{\beta}\|_2 \leq \left(\frac{32\rho_{\hat{A},k+\ell}^{(2)}}{t\mu_{\hat{A},k+\ell}^{(2)}} + 4\right)\frac{r_k^{(1)}(\bar{\beta})}{\sqrt{\ell}} + \left(\frac{8\rho_{\hat{A},k+\ell}^{(2)}}{t\mu_{\hat{A},k+\ell}^{(2)}} + 4\right)r_k^{(2)}(\bar{\beta}) + \frac{8}{t\mu_{\hat{A},k+\ell}^{(2)}}\sqrt{k}\lambda.$$

Corollary 6.1 is directly comparable to the main result of [4] for the Dantzig selector, which is given by the estimator

$$\hat{\beta}_D = \underset{\beta \in R^d}{\arg\min} \|\beta\|_1 \qquad \text{subject to } \sup_j \left|\sum_{i=1}^n \mathbf{x}_{i,j}(\mathbf{x}_i^T\beta - \mathbf{y}_i)\right| \leq b_D.$$

Their main result is stated below, in Theorem 6.1. It uses a different quantity $\bar{\theta}_{A,k,\ell}$, which is defined as

$$\bar{\theta}_{A,k,\ell} = \sup_{\beta \in R^\ell, I, J} \frac{\|A_{I,J}\beta\|_2}{\|\beta\|_2},$$

using notation of Definition 3.1. It is easy to see that $\theta_{A,k,\ell}^{(2)} \leq \bar{\theta}_{A,k,\ell}\sqrt{\ell}$.

THEOREM 6.1 [4]. *Assume that there exists a vector* $\bar{\beta} \in R^d$ *with* $s$ *nonzero components, such that* $\mathbf{y}_i = \bar{\beta}^T \mathbf{x}_i + \varepsilon_i$, *where* $\varepsilon_i \sim N(0,\sigma^2)$ *are i.i.d.*

*Gaussian noises. Let $\hat{A} = \frac{1}{n}\sum_{i=1}^{n}\mathbf{x}_i\mathbf{x}_i^T$, and assume that $\hat{A}_{j,j} \leq 1$ for all $j$. Given $t_D > 0$ and $\delta \in (0,1)$, we set $b_D = \sqrt{n}\lambda_D\sigma$, with*

$$\lambda_D = (\sqrt{1 - (\ln\delta + \ln(\sqrt{\pi\ln d}))/\ln d} + t_D^{-1})\sqrt{2\ln d}.$$

*Let $\bar{\theta}_{\hat{A},2s} = \max(\rho_{\hat{A},2s}^{(2)} - 1, 1 - \mu_{\hat{A},2s}^{(2)})$. If $\bar{\theta}_{\hat{A},2s} + \bar{\theta}_{\hat{A},2s,s} < 1 - t_D$, then, with probability exceeding $1 - \delta$,*

$$\|\hat{\beta}_D - \bar{\beta}\|_2^2 \leq C_2(\bar{\theta}_{\hat{A},2s}, \bar{\theta}_{\hat{A},2s,s})\lambda_D^2\left(\frac{(k+1)\sigma^2}{n} + r_k^{(2)}(\bar{\beta})^2\right).$$

*The quantity $C_2(a,b)$ is defined as $C_2(a,b) = \frac{2C_0(a,b)}{1-a-b} + \frac{2b(1+a)}{(1-a-b)^2} + \frac{1+a}{1-a-b}$, where $C_0(a,b) = 2\sqrt{2}(1 + \frac{1-a^2}{1-a-b}) + (1 + 1/\sqrt{2})(1+a)^2/(1-a-b).$*

In order to see that Corollary 6.1 is comparable to Theorem 6.1, we shall compare their conditions and consequences. To this end, we can pick any $\ell \in [k, s]$ in Corollary 6.1.

We shall first look at the conditions. The condition required in Theorem 6.1 is $\bar{\theta}_{\hat{A},2s} + \bar{\theta}_{\hat{A},2s,s} < 1 - t_D$, which implies that $\bar{\theta}_{\hat{A},2s,s} < \mu_{\hat{A},2s}^{(2)} - t_D$. This condition is stronger than $\theta_{\hat{A},k+\ell,\ell}^{(2)}/\sqrt{k} \leq \bar{\theta}_{\hat{A},k+\ell,\ell} < \mu_{\hat{A},k+\ell}^{(2)} - t_D$, which implies the condition $t > 0$ in Corollary 6.1:

$$t = 1 - \pi_{\hat{A},k+\ell,\ell}^{(2)}\frac{\sqrt{k}}{\ell} \geq 1 - \frac{\theta_{\hat{A},k+\ell,\ell}^{(2)}\sqrt{k}}{\mu_{\hat{A},k+\ell}^{(2)}\ell} \geq 1 - \frac{(\mu_{\hat{A},k+\ell}^{(2)} - t_D)k}{\mu_{\hat{A},k+\ell}^{(2)}\ell}$$

$$= (1 - k/\ell) + \frac{t_D k}{\mu_{\hat{A},k+\ell}^{(2)}\ell} > 0.$$

Therefore, Corollary 6.1 can be applied, as long as Theorem 6.1 can be applied. Moreover, as long as $t_D > 0$, $t$ is never much smaller than $t_D$ but can be significantly larger (e.g., $t > 0.5$ when $k \leq \ell/2$, even if $t_D$ is close to zero). It is also obvious that the condition $t > 0$ does not imply that $t_D > 0$. Therefore, the condition of Corollary 6.1 is strictly weaker. As discussed in Section 4, if $\rho_{A,\ell}^{(2)}/\mu_{A,k+2\ell}^{(2)} \leq c$ for a constant $c > 0$ when $\ell \leq ck$; then, the condition $t > 0$ holds with $\ell = ck$.

Next, we shall look at the consequences of the two theorems when both $t > 0$ and $t_D > 0$. Ignoring constants, the bound in Theorem 6.1, with $\lambda_D = O(\sqrt{\ln(d/\delta)})$, can be written as

$$\|\hat{\beta}_D - \bar{\beta}\|_2 = O(\sqrt{\ln(d/\delta)})(r_k^{(2)}(\bar{\beta}) + \sigma\sqrt{k/n}).$$

However, the proof itself implies a stronger bound of the form

$$(3) \qquad \|\hat{\beta}_D - \bar{\beta}\|_2 = O(\sigma\sqrt{k\ln(d/\delta)/n} + r_k^{(2)}(\bar{\beta})).$$

In comparison, in Corollary 6.1, we can pick $\lambda = O(\sigma\sqrt{\ln(d/\delta)/n})$, and then the bound can be written as (with $\ell = s$)

$$(4) \qquad \|\hat{\beta} - \bar{\beta}\|_2 = O(\sigma\sqrt{k\ln(d/\delta)/n} + r_k^{(2)}(\bar{\beta}) + r_k^{(1)}(\bar{\beta})/\sqrt{s}).$$

Note that we do not have to assume that $\bar{\beta}$ only contains $s$ nonzero components. The quantity $r_k^{(1)}(\bar{\beta})/\sqrt{s}$ is no more than $r_k^{(2)}(\bar{\beta})$ under the sparsity of $\bar{\beta}$ assumed in Theorem 6.1. It is thus clear that (4) has a more general form than that of (3).

It was pointed out in [1] that Lasso and the Dantzig selector are quite similar, and the authors presented a simultaneous analysis of both. Since the explicit parameter estimation bounds in [1] are with the case $k = \mathrm{supp}_0(\bar{\beta})$, it is natural to ask whether our results (in particular, the first claim of Theorem 4.1) can also be applied to the Dantzig selector, so that a simultaneous analysis similar to that of [1] can be established. Unfortunately, the techniques used in this paper do not immediately give an affirmative answer. This is because a Lasso-specific property is used in our proof of Lemma 10.4, and the property does not hold for the Dantzig selector. However, we conjecture that it may still be possible to prove similar results for the Dantzig selector through different techniques such as those employed in [4] and [6].

**7. Feature selection through coefficient thresholding.** A fundamental result of $L_1$ regularization is its feature selection consistency property, which is considered in [16] and more formally analyzed in [20]. It was shown that, under a strong *irrepresentable condition* (introduced in [20]) together with the sparse eigenvalue condition, the set $\mathrm{supp}_0(\hat{\beta})$, with $\hat{\beta}$ estimated using Lasso, may be able to consistently identify features with coefficients larger than a threshold of order $\sqrt{k}\lambda$ (with $\lambda = O(\sigma\sqrt{\ln(d/\delta)/n})$). Here, $k$ is the sparsity of the true target. That is, with probability $1 - \delta$, all coefficients larger than a certain threshold of the order $O(\sigma\sqrt{k\ln(d/\delta)/n})$ remain nonzero, while all zero coefficients remain zero. It was also shown that a slightly weaker irrepresentable condition is necessary for Lasso to possess this property. For Lasso, the $\sqrt{k}$ factor cannot be removed (under the sparse eigenvalue assumption plus the irrepresentable condition) unless additional conditions (such as the mutual incoherence assumption in Corollary 7.1) are imposed. Also, see [5, 18] for related results without the $\sqrt{k}$ factor.

It was acknowledged in [19] and [17] that the irrepresentable condition can be quite strong (e.g., often more restricted than eigenvalue conditions required for Corollary 6.1). This is the motivation of the sparse eigenvalue condition introduced in [19], although such a condition does not necessarily yield consistent feature selection under the scheme of [20], which employs the set $\mathrm{supp}_0(\hat{\beta})$ to identify features. However, limitations of the irrepresentable condition can be removed by considering $\mathrm{supp}_\alpha(\hat{\beta})$ with $\alpha > 0$.

In this section, we consider a more extended view of feature selection, where a practitioner would like to find relevant features with coefficient magnitude larger than some threshold $\alpha$ that is not necessarily zero. Features with small coefficients are regarded as irrelevant features, which are not distinguished from zero for practical purposes. The threshold $\alpha$ can be pre-chosen based on the interests of the practitioner as well as our knowledge of the underlying problem. We are interested in the relationship of features estimated from the solution $\hat{\beta}$ of (1) and the true relevant features obtained from $\bar{\beta}$. The following result is a simple consequence of Theorem 4.1, where we use $\|\hat{\beta} - \bar{\beta}\|_p$ for some large $p$ to approximate $\|\hat{\beta} - \bar{\beta}\|_\infty$ (which is needed for feature selection). A consequence of the result (see Corollary 7.2) is that, using a nonzero threshold $\alpha$ (rather than zero-threshold of [20]), it is possible to achieve consistent feature selection even if the irrepresentable condition in [20] is violated. For clarity, we choose a simplified statement with sparse target $\bar{\beta}$. However, it is easy to see from the proof that, just as in Theorem 4.1, a similar but more complicated statement holds, even when the target is not sparse.

THEOREM 7.1. *Let Assumption 4.1 hold, and let $\hat{A} = \frac{1}{n}\sum_{i=1}^n \mathbf{x}_i \mathbf{x}_i^T$ and $a = (\sup_j \hat{A}_{j,j})^{1/2}$. Let $\bar{\beta} \in R^d$ be the true target vector with $\mathbf{E}y = \bar{\beta}^T \mathbf{x}$, and assume that $|\operatorname{supp}_0(\bar{\beta})| = k$. Let $\hat{\beta}$ be the solution of (1). Given $\delta \in (0, 1)$, then, with probability larger than $1 - \delta$, the following claim is true. For all $\varepsilon \in (0, 1)$, if there exist $(d - k)/2 \geq \ell \geq k$, $t \in (0, 1)$, and $p \in [1, \infty]$ so that:*

- $\lambda \geq 4(2 - t)t^{-1}(\sigma a\sqrt{2\ln(2d/\delta)/n})$;
- *either $8(\varepsilon\alpha\omega_{\hat{A},k+\ell}^{(p)})^{-1}k^{1/p}\lambda \leq t \leq 1 - \pi_{\hat{A},k+\ell,\ell}^{(p)}k^{1-1/p}/\ell$, or $8(\varepsilon\alpha\mu_{\hat{A},k+\ell}^{(p)})^{-1}(k + \ell)^{1/p}\lambda \leq t \leq 1 - \gamma_{\hat{A},k+\ell,\ell}^{(p)}k^{1-1/p}/\ell$;*

*then $\operatorname{supp}_{(1+\varepsilon)\alpha}(\bar{\beta}) \subset \operatorname{supp}_\alpha(\hat{\beta}) \subset \operatorname{supp}_{(1-\varepsilon)\alpha}(\bar{\beta})$.*

If either $t \leq 1 - \pi_{\hat{A},k+\ell,\ell}^{(p)}k^{1-1/p}/\ell$ or $t \leq 1 - \gamma_{\hat{A},k+\ell,\ell}^{(p)}k^{1-1/p}/\ell$, then the result can be applied as long as $\alpha$ is sufficiently large. As we have pointed out after Theorem 4.1, if the sparse eigenvalue condition holds at block size of order $k^{2-2/p}$ for some $p \geq 2$, then one can take $\ell = O(k^{2-2/p})$, so that the condition $t \leq 1 - \gamma_{\hat{A},k+\ell,\ell}^{(p)}k^{1-1/p}/\ell$ is satisfied. This implies that we may take $\alpha = O(k^{2/p-2/p^2}\lambda) = O(\sigma k^{2/p-2/p^2}\sqrt{\ln(d/\delta)/n})$, assuming that $\mu_{\hat{A},k+\ell}^{(p)}$ is bounded from below (which holds when $\hat{A}$ is $p$-norm diagonal dominant at size $k+\ell$, according to Proposition 3.1). That is, sparse eigenvalue condition at a certain block size of order $k^{2-2/p}$, together with the boundedness of $\mu_{\hat{A},k+\ell}^{(p)}$, imply that one can distinguish coefficients of magnitude larger than a threshold of order $\sigma k^{2/p-2/p^2}\sqrt{\ln(d/\delta)/n}$ from zero. In particular, if $p=\infty$, we can distinguish nonzero coefficients of order $\sigma\sqrt{\ln(d/\delta)/n}$ from zero. For simplicity, we state such a result under the mutual incoherence assumption.

COROLLARY 7.1.  *Let Assumption 4.1 hold, let $\hat{A} = \frac{1}{n}\sum_{i=1}^{n}\mathbf{x}_i\mathbf{x}_i^T$ and assume that $\hat{A}_{j,j} = 1$ for all $j$. Define $M_{\hat{A}} = \sup_{i\neq j}|\hat{A}_{i,j}|$. Let $\hat{\beta}$ be the solution of (1). Let $\bar{\beta} \in R^d$ be the true target vector with $\mathbf{E}y = \bar{\beta}^T\mathbf{x}$ and $k = |\operatorname{supp}_0(\bar{\beta})|$. Assume that $kM_{\hat{A}} \leq 0.25$ and $3k \leq d$. Given $\delta \in (0,1)$, with probability larger than $1-\delta$, if $\alpha/32 \geq \lambda \geq 12\sigma\sqrt{2\ln(2d/\delta)/n}$, then $\operatorname{supp}_{(1+\varepsilon)\alpha}(\bar{\beta}) \subset \operatorname{supp}_{\alpha}(\hat{\beta}) \subset \operatorname{supp}_{(1-\varepsilon)\alpha}(\bar{\beta})$, where $\varepsilon = 32\lambda/\alpha$.*

One can also obtain a formal result on the asymptotic consistency of feature selection. An example is given below. In the description, we allow the problem to vary with sample size $n$, and study the asymptotic behavior when $n \to \infty$. Therefore, except for the input vectors $\mathbf{x}_i$, all other quantities such as $d$, $\bar{\beta}$, etc., will be denoted with subscript $n$. The input vectors $\mathbf{x}_i \in R^{d_n}$ also vary with $n$; however, we drop the subscript $n$ to simplify the notation. The statement of our result is in the same style as a corresponding asymptotic feature selection consistency theorem of [20] for the zero-thresholding scheme $\operatorname{supp}_0(\hat{\beta})$, which requires the stronger irrepresentable condition in addition to the sparse eigenvalue condition. In contrast, our result employs nonzero thresholding $\operatorname{supp}_{\alpha_n}(\hat{\beta})$, with an appropriately chosen sequence of decreasing $\alpha_n$; the result only requires the sparse eigenvalue condition (and, for clarity, we only consider $p=2$ instead of general $p$ discussed above) without the need for irrepresentable condition.

COROLLARY 7.2.  *Consider regression problems indexed by the sample size $n$, and let the corresponding true target vector be $\bar{\beta}_n = [\bar{\beta}_{n,1},\ldots,\bar{\beta}_{n,d_n}] \in R^{d_n}$, where $\mathbf{E}y = \bar{\beta}_n^T\mathbf{x}$. Let Assumption 4.1 hold, with $\sigma$ independent of $n$. Assume that there exists $a > 0$ that is independent of $n$, such that $\frac{1}{n}\sum_{i=1}^{n}\mathbf{x}_{i,j}^2 \leq a^2$ for all $j$. Denote, by $\hat{\beta}_n$, the solution $\hat{\beta}$ of (1) with $\lambda = 12\sigma a\sqrt{2(\ln(2d_n) + n^{s'})/n}$, where $s' \in (0,1)$. Pick $s \in (0, 1-s')$, and set $\alpha_n = n^{-s/2}$. Then, as $n \to \infty$, $P(\operatorname{supp}_{\alpha_n}(\hat{\beta}_n) \neq \operatorname{supp}_0(\bar{\beta}_n)) = O(\exp(-n^{s'}))$ if the following conditions hold:*

1. *$\bar{\beta}_n$ only has $k_n = o(n^{1-s-s'})$ nonzero coefficients;*
2. *$k_n\ln(d_n) = o(n^{1-s})$;*
3. *$1/\min_{j\in\operatorname{supp}_0(\bar{\beta}_n)}|\bar{\beta}_{n,j}| = o(n^{s/2})$;*

4. Let $\hat{A}_n = \frac{1}{n}\sum_{i=1}^n \mathbf{x}_i\mathbf{x}_i^T \in R^{d_n \times d_n}$. There exists a positive integer $q_n$ such that $(1 + 2q_n)k_n \leq d_n$, $1/\mu_{\hat{A}_n,(1+2q_n)k_n}^{(2)} = O(1)$, and $\rho_{\hat{A}_n,(1+2q_n)k_n}^{(2)} \leq (1 + q_n)\mu_{\hat{A}_n,(1+2q_n)k_n}^{(2)}$.

The conditions of the corollary are all standard. Similar conditions have also appeared in [20]. The first condition simply requires $\bar{\beta}_n$ to be sufficiently sparse; if $k_n$ is in the same order of $n$, then one cannot obtain meaningful consistency results. The second condition requires that $d_n$ is not too large, and, in particular, that it should be sub-exponential in $n$; otherwise, our analysis does not lead to consistency. The third condition requires that $|\bar{\beta}_{n,j}|$ be sufficiently large when $j \in \mathrm{supp}_0(\bar{\beta}_n)$. In particular, the condition implies that each feature component $|\bar{\beta}_{n,j}|$ needs to be larger than the 2-norm noise level $\sigma\sqrt{k_n \ln(d_n)/n}$. If some component $\bar{\beta}_{n,j}$ is too small, then we cannot distinguish it from the noise. Note that, since the 2-norm parameter estimation bound is used here, we have a $\sqrt{k_n}$ factor in the noise level. Under stronger conditions, such as mutual incoherence, this $\sqrt{k_n}$ factor can be removed (as shown in Corollary 7.1). Finally, the fourth condition is the sparse eigenvalue assumption; it can also be replaced by some other conditions (such as mutual incoherence). In comparison, [20] employed zero-threshold scheme with $\alpha_n = 0$; therefore, in addition to our assumptions, the irrepresentable condition is also required.

## 8. Two-stage $L_1$ regularization with selective penalization.
We shall refer to the feature components corresponding to the large coefficients as *relevant features* and the feature components smaller than an appropriately defined cut-off threshold $\alpha$ as *irrelevant features*. Theorem 7.1 implies that Lasso can be used to approximately identify the set of relevant features $\mathrm{supp}_\alpha(\bar{\beta})$. This property can be used to improve the standard Lasso. In this context, we observe that as an estimation method, $L_1$ regularization has two important properties, which are as follows:

1. Shrink estimated coefficients corresponding to irrelevant features toward zero;
2. Shrink estimated coefficients corresponding to relevant features toward zero.

While the first effect is desirable, the second effect is not. In fact, we should avoid shrinking the coefficients corresponding to the relevant features if we can identify these features. In this case, the standard $L_1$ regularization may have sub-optimal performance. In order to improve $L_1$, we observe that, under appropriate conditions such as those of Theorem 7.1, estimated coefficients corresponding to relevant features tend to be larger than estimated coefficients corresponding to irrelevant features. Therefore, after the first stage of $L_1$ regularization, we can identify the relevant features by picking the components corresponding to the largest coefficients. Those coefficients are over-shrinked in the first stage. This problem can be fixed by applying a second stage of $L_1$ regularization, where we do not penalize the features selected in the first stage. The procedure is described in Figure 1. Its overall effect is to "unshrink" coefficients of relevant features identified in the first stage. In practice, instead of tuning $\alpha$, we may also let $\mathrm{supp}_\alpha(\hat{\beta})$ contain exactly $q$ elements, and simply tune the integer valued $q$. The parameters can then be tuned by cross-validation in sequential order: first, find $\lambda$ to optimize stage 1 prediction accuracy; second, find $q$ to optimize stage 2 prediction accuracy. If cross-validation works well, then this tuning method ensures that the two-stage selective penalization procedure is never much worse than the one-stage procedure in practice, because they are equivalent with $q = 0$. However, under the right conditions, we can prove a much better bound for this two stage procedure, as shown in Theorem 8.1.

A related method, called *relaxed Lasso*, was proposed recently by Meinshausen [15], which is similar to a two-stage Dantzig selector in [4] (also see [12] for a more detailed study). Their idea differs from our proposal in that, in the second stage, the parameter coefficients $\beta'_j$ are forced to be zero when $j \notin \mathrm{supp}_0(\hat{\beta})$. It was pointed out in [15] that, if $\mathrm{supp}_0(\hat{\beta})$ can exactly identify all nonzero components of the target vector, then, in the second stage, the relaxed Lasso can asymptotically remove the bias in the first-stage Lasso. However, it is not clear what theoretical result can be stated when Lasso cannot exactly identify all relevant features. In the general case, it is not easy to ensure that relaxed Lasso does not degrade the performance when some relevant coefficients become zero in the first stage. On the contrary, the two-stage selective penalization procedure in Figure 1 does not require that all relevant features are identified. Consequently, we are able to prove a result for Figure 1 with no counterpart for relaxed Lasso. For clarity, the

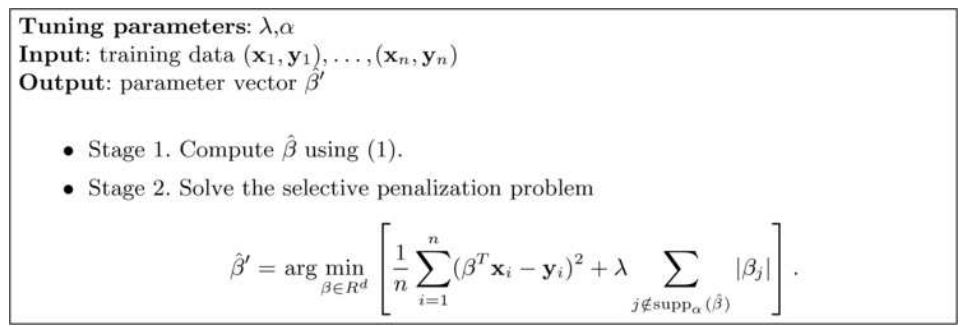

**Tuning parameters**: $\lambda, \alpha$
**Input**: training data $(\mathbf{x}_1, \mathbf{y}_1), \ldots, (\mathbf{x}_n, \mathbf{y}_n)$
**Output**: parameter vector $\hat{\beta}'$

- Stage 1. Compute $\hat{\beta}$ using (1).
- Stage 2. Solve the selective penalization problem

$$\hat{\beta}' = \arg\min_{\beta \in R^d} \left[ \frac{1}{n} \sum_{i=1}^{n} (\beta^T \mathbf{x}_i - \mathbf{y}_i)^2 + \lambda \sum_{j \notin \mathrm{supp}_\alpha(\hat{\beta})} |\beta_j| \right].$$

FIG. 1. *Two-stage $L_1$ regularization with selective penalization.*

result is stated under similar conditions to those of the Dantzig selector in Theorem 6.1 with sparse targets and $p = 2$ only. Both restrictions can be easily removed, with a more complicated version of Theorem 7.1, to deal with nonsparse targets (which can be easily obtained from Theorem 4.1) as well as the general form of Lemma 10.4, which allows $p \neq 2$.

THEOREM 8.1. *Let Assumption 4.1 hold, and let* $\hat{A} = \frac{1}{n} \sum_{i=1}^{n} \mathbf{x}_i \mathbf{x}_i^T$ *and* $a = (\sup_j \hat{A}_{j,j})^{1/2}$. *Consider any target vector* $\bar{\beta} \in R^d$ *such that* $\mathbf{E}y = \bar{\beta}^T \mathbf{x}$, *and* $\bar{\beta}$ *contains only s nonzeros. Let* $k = |\mathrm{supp}_\lambda(\bar{\beta})|$. *Consider the two-stage selective penalization procedure in Figure 1. Given* $\delta \in (0, 0.5)$, *with probability larger than* $1 - 2\delta$ *for all* $(d-s)/2 \geq \ell \geq s$ *and* $t \in (0,1)$, *assume the following:*

- $t \leq 1 - \pi_{\hat{A}, k+\ell, \ell}^{(2)} k^{0.5}/\ell$;
- $0.5\alpha \geq \lambda \geq 4(2-t) t^{-1} \sigma a \sqrt{2\ln(2d/\delta)/n}$;
- *either* $16(\alpha \omega_{\hat{A}, s+\ell}^{(p)})^{-1} s^{1/p} \lambda \leq t \leq 1 - \pi_{\hat{A}, s+\ell, \ell}^{(p)} s^{1-1/p}/\ell$, *or* $16(\alpha \mu_{\hat{A}, s+\ell}^{(p)})^{-1} (s + \ell)^{1/p} \lambda \leq t \leq 1 - \gamma_{\hat{A}, s+\ell, \ell}^{(p)} s^{1-1/p}/\ell$.

*Then,*

$$\|\hat{\beta}' - \bar{\beta}\|_2 \leq \frac{8}{t\mu_{\hat{A}, k+\ell}^{(2)}} [5\rho_{\hat{A}, k+\ell}^{(2)} r_k^{(2)}(\bar{\beta}) + \sqrt{k-q}\lambda + a\sigma(1 + \sqrt{20\ln(1/\delta)})\sqrt{q/n}]$$
$$+ 8r_k^{(2)}(\bar{\beta}),$$

*where* $q = |\mathrm{supp}_{1.5\alpha}(\bar{\beta})|$.

Again, we include a simplification of Theorem 8.1 under the mutual incoherence condition.

COROLLARY 8.1. *Let Assumption 4.1 hold, let* $\hat{A} = \frac{1}{n} \sum_{i=1}^{n} \mathbf{x}_i \mathbf{x}_i^T$ *and assume that* $\hat{A}_{j,j} = 1$ *for all j. Define* $M_{\hat{A}} = \sup_{i \neq j} |\hat{A}_{i,j}|$. *Consider any target vector* $\bar{\beta}$ *such that* $\mathbf{E}y = \bar{\beta}^T \mathbf{x}$, *and assume that* $\bar{\beta}$ *contains only s nonzeros where* $s \leq d/3$ *and assume that* $M_{\hat{A}} s \leq 1/6$. *Let* $k = |\mathrm{supp}_\lambda(\bar{\beta})|$. *Consider the two-stage selective penalization procedure in Figure 1. Given* $\delta \in (0, 0.5)$, *with probability larger than* $1 - 2\delta$, *if* $\alpha/48 \geq \lambda \geq 12\sigma\sqrt{2\ln(2d/\delta)/n}$, *then*

$$\|\hat{\beta}' - \bar{\beta}\|_2 \leq 24\sqrt{k-q}\lambda + 24\sigma\left(1 + \sqrt{\frac{20q}{n}\ln(1/\delta)}\right) + 168\delta_k^{(2)}(\bar{\beta}),$$

*where* $q = |\mathrm{supp}_{1.5\alpha}(\bar{\beta})|$.

Theorem 8.1 can significantly improve the corresponding one-stage result (see Corollary 6.1 and Theorem 6.1) when $r_k^{(2)}(\bar{\beta}) \ll \sqrt{k}\lambda$ and $k - q \ll k$. The latter condition is true when $|\operatorname{supp}_{1.5\alpha}(\bar\beta)| \approx |\operatorname{supp}_\lambda(\bar\beta)|$. In such a case, we can identify most features in $\operatorname{supp}_\lambda(\bar\beta)$. These conditions are satisfied when most nonzero coefficients in $\operatorname{supp}_\lambda(\bar\beta)$ are relatively large in magnitude and the other coefficients are small (in 2-norm). That is, the two-stage procedure is superior when the target $\bar\beta$ can be decomposed as a sparse vector with large coefficients plus another (less sparse) vector with small coefficients. In the extreme case, when $r_k^{(2)}(\beta) = 0$ and $q = k$, we obtain $\|\hat\beta' - \bar\beta\|_2 = O(\sqrt{k\ln(1/\delta)/n})$ instead of $\|\hat\beta - \bar\beta\|_2 = O(\sqrt{k\ln(d/\delta)/n})$ for the one-stage Lasso. The difference can be significant when $d$ is large.

Finally, we shall point out that the two-stage selective penalization procedure may be regarded as a two-step approximation to solving the least squares problem with a nonconvex regularization:

$$\hat\beta' = \underset{\beta \in R^d}{\arg\min} \left[ \frac{1}{n} \sum_{i=1}^{n} (\beta^T \mathbf{x}_i - \mathbf{y}_i)^2 + \lambda \sum_{j=1}^{d} \min(\alpha, |\beta_j|) \right].$$

However, for high-dimensional problems, it not clear whether one can effectively find a good solution using such a nonconvex regularization condition. When $d$ is sufficiently large, one can often find a vector $\beta$, such that $|\beta_j| > \alpha$ and it perfectly fits (thus overfits) the data. This $\beta$ is clearly a local minimum for this nonconvex regularization condition, since the regularization has no effect locally for such a vector $\beta$. Therefore, the two-stage $L_1$ approximation procedure in Figure 1 not only preserves desirable properties of convex programming, but also prevents such a local minimum to contaminate the final solution.

**9. Experiments.** Although our investigation is mainly theoretical, it is useful to verify whether the two stage procedure can improve the standard Lasso in practice. In the following, we show with a synthetic data and a real data that the two-stage procedure can be helpful. Although more comprehensive experiments are still required, these simple experiments show that the two-stage method is useful at least on datasets with the right properties, which is consistent with our theory. Note that, instead of tuning the $\alpha$ parameter in Figure 1, in the following experiments, we tune the parameter $q = \operatorname{supp}_\alpha(\hat\beta)$, which is more convenient. The standard Lasso corresponds to $q = 0$.

9.1. *Simulation data.* In this experiment, we generate an $n \times d$ random matrix with its column $j$ corresponding to $[\mathbf{x}_{1,j}, \ldots, \mathbf{x}_{n,j}]$, and each element of the matrix is an independent standard Gaussian $N(0, 1)$. We then normalize its columns so that $\sum_{i=1}^{n} \mathbf{x}_{i,j}^2 = n$. A truly sparse target $\bar\beta$, is generated with $k$ nonzero elements that are uniformly distributed from $[-10, 10]$. Observe that $\mathbf{y}_i = \bar\beta^T \mathbf{x}_i + \varepsilon_i$, where each $\varepsilon_i \sim N(0, \sigma^2)$. In this experiment, we take $n = 25, d = 100, k = 5$ and $\sigma = 1$, and repeat the experiment 100 times. The average training error and parameter estimation error in 2-norm are reported in Figure 2. We compare the performance of the two-stage method with different $q$ versus the regularization parameter $\lambda$. Clearly, the training error becomes smaller when $q$ increases. The smallest estimation error for this example is achieved at $q = 3$. This shows that the two-stage procedure with appropriately chosen $q$ performs better than the standard Lasso (which corresponds to $q = 0$).

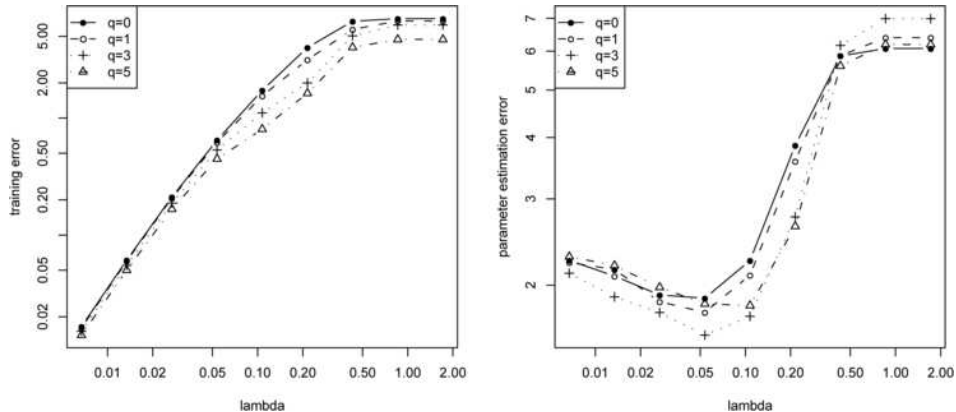

Fig. 2. *Performance of the algorithms on simulation data. Left: average training squared error versus $\lambda$; right: parameter estimation error versus $\lambda$.*

9.2. *Real data.* We use real data to illustrate the effectiveness of two-stage $L_1$ regularization. For simplicity, we only report the performance on a single data, *Boston Housing*. This is the housing data for 506 census tracts in Boston from the 1970 census, available from the *UCI machine learning database repository* (http://archive.ics.uci.edu/ml/). Each census tract is a data-point with 13 features (we add a constant offset on e as the 14th feature), and the desired output is the housing price. In the experiment, we randomly partition the data into 20 training plus 456 test points. We perform the experiments 100 times and report training and test squared error versus the regularization parameter $\lambda$ for different $q$. The results are plotted in Figure 3. In this case, $q = 1$ achieves the best performance. Note that this dataset contains only a small number ($d = 14$) features, which is not the case we are interested in (most of other UCI data similarly contain only small number of features). In order to illustrate the advantage of the two-stage method more clearly, we also consider a modified Boston Housing data, where we append 20 random features (similar to the simulation experiments) to the original Boston Housing data, and rerun the experiments. The results are shown in Figure 4. As we can expect, the effect of using $q > 0$ becomes far more apparent. This again verifies that the two-stage method can be superior to the standard Lasso ($q = 0$) on some data.

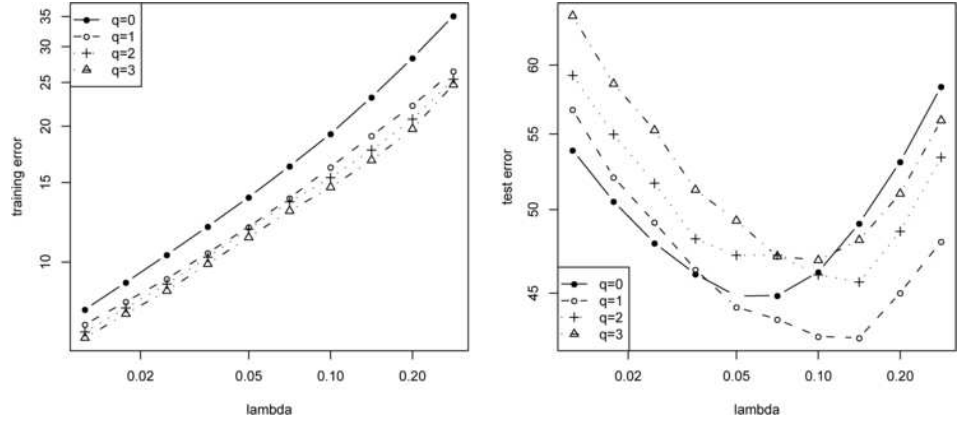

Fig. 3.   *Performance of the algorithms on the original Boston Housing data. Left: average training squared error versus λ; right: test squared error versus λ.*

**10. Proofs.**   In the proof, we use the following convention: let $I$ be a subset of $\{1, \ldots, d\}$ and a vector $\beta \in R^d$, then $\beta_I$ denotes either the restriction of $\beta$ to indices $I$, which lies in $R^{|I|}$, or its embedding into the original space $R^d$ with components not in $I$ set to zero.

10.1. *Proof of Proposition 3.1.*   Given any $\mathbf{v} \in R^k$ and $\mathbf{u} \in R^\ell$, without loss of generality, we may assume that $\|\mathbf{v}\|_2 = 1$ and $\|\mathbf{u}\|_2 = 1$ in the following derivation. Take indices $I$ and $J$ as in the definition. We let $I' = I \cup J$. Given any $\alpha \in R$, let $\mathbf{u}'^T = [\mathbf{v}^T, \alpha \mathbf{u}^T] \in R^{\ell+k}$. By definition, we have

$$\mu_{A,\ell+k}^{(2)} \|\mathbf{u}'\|_2^2 \leq \mathbf{u}'^T A_{I',I'} \mathbf{u}' = \mathbf{v}^T A_{I,I} \mathbf{v} + 2\alpha \mathbf{v}^T A_{I,J} \mathbf{u} + \alpha^2 \mathbf{u}^T A_{J,J} \mathbf{u}.$$

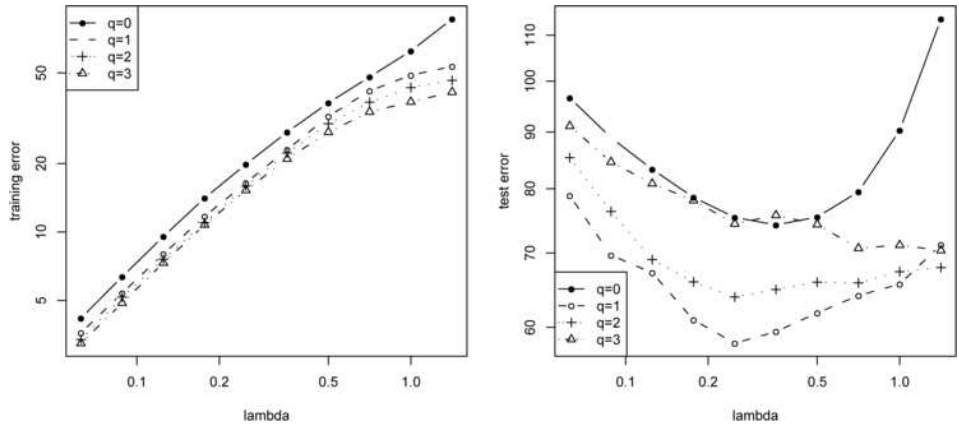

Fig. 4.   *Performance of the algorithms on the modified Boston Housing data. Left: average training squared error versus λ; right: test squared error versus λ.*

Let $b = \mathbf{v}^T A_{I,J} \mathbf{u}$, $c_1 = \mathbf{v}^T A_{I,I} \mathbf{v}$ and $c_2 = \mathbf{u}^T A_{J,J} \mathbf{u}$. The above inequality can be written as

$$\mu_{A,\ell+k}^{(2)}(1 + \alpha^2) \leq c_1 + 2\alpha b + \alpha^2 c_2.$$

By optimizing over $\alpha$, we obtain $(c_1 - \mu_{A,\ell+k}^{(2)})(c_2 - \mu_{A,\ell+k}^{(2)}) \geq b^2$. Therefore, $b^2 \leq (\rho_{A,\ell}^{(2)} - \mu_{A,\ell+k}^{(2)})(\rho_{A,k}^{(2)} - \mu_{A,\ell+k}^{(2)})$, which implies that (with $\|\mathbf{v}\|_2 = \|\mathbf{u}\|_2 = 1$)

$$\frac{\mathbf{v}^T A_{I,J} \mathbf{u}}{\|\mathbf{v}\|_2 \|\mathbf{u}\|_\infty} \leq \frac{|\mathbf{v}^T A_{I,J} \mathbf{u}|}{\|\mathbf{v}\|_2 \|\mathbf{u}\|_2 / \sqrt{\ell}} \leq \sqrt{\ell} |b|$$

$$\leq \sqrt{\ell} \sqrt{(\rho_{A,\ell}^{(2)} - \mu_{A,\ell+k}^{(2)})(\rho_{A,k}^{(2)} - \mu_{A,\ell+k}^{(2)})}.$$

Since $\mathbf{v}$ and $\mathbf{u}$ are arbitrary, this implies the first inequality.

The second inequality follows from $\|\mathbf{v}\|_p \leq k^{\max(0,1/p-0.5)} \|\mathbf{v}\|_2$ for all $\mathbf{v} \in R^k$, so that

$$\frac{\|A_{I,J} \mathbf{u}\|_p}{\|\mathbf{u}\|_\infty} \leq k^{\max(0,1/p-0.5)} \frac{\|A_{I,J} \mathbf{u}\|_2}{\|\mathbf{u}\|_\infty}.$$

From $(c_1 - \mu_{A,\ell+k}^{(2)})(c_2 - \mu_{A,\ell+k}^{(2)}) \geq b^2$, we also obtain

$$4b^2/c_1^2 \leq 4c_1^{-1}(1 - \mu_{A,\ell+k}^{(2)}/c_1)(c_2 - \mu_{A,\ell+k}^{(2)}) \leq (c_2 - \mu_{A,\ell+k}^{(2)})/\mu_{A,\ell+k}^{(2)}$$

$$\leq \rho_{A,\ell}^{(2)}/\mu_{A,\ell+k}^{(2)} - 1.$$

Note that, in the above derivation, we have used $4\mu_{A,\ell+k}^{(2)} c_1^{-1}(1 - \mu_{A,\ell+k}^{(2)} c_1^{-1}) \leq 1$. Therefore, with $\|\mathbf{v}\|_2 = \|\mathbf{u}\|_2 = 1$,

$$\frac{\mathbf{v}^T A_{I,J} \mathbf{u} \|\mathbf{v}\|_2}{\mathbf{v}^T A_{I,I} \mathbf{v} \|\mathbf{u}\|_\infty} \leq \frac{|\mathbf{v}^T A_{I,J} \mathbf{u}|}{\mathbf{v}^T A_{I,I} \mathbf{v}/\sqrt{\ell}} = \frac{|b|}{c_1} \sqrt{\ell} \leq 0.5\ell^{1/2} \sqrt{\rho_{A,\ell}^{(2)}/\mu_{A,\ell+k}^{(2)} - 1}.$$

Because $\mathbf{v}$ and $\mathbf{u}$ are arbitrary, we obtain the third inequality.

The fourth inequality follows from $\max(0, \mathbf{v}^T A_{I,I} \mathbf{v}^{p-1}) \geq \omega_{A,k}^{(p)} \|\mathbf{v}\|_p^p$, so that

$$\frac{(\mathbf{v}^{p-1})^T A_{I,J} \mathbf{u} \|\mathbf{v}\|_p}{\max(0, \mathbf{v}^T A_{I,I} \mathbf{v}^{p-1}) \|\mathbf{u}\|_\infty} \leq \frac{1}{\omega_{A,k}^{(p)}} \frac{|(\mathbf{v}^{p-1})^T A_{I,J} \mathbf{u}|}{\|\mathbf{v}\|_p^{p-1} \|\mathbf{u}\|_\infty} \leq \frac{1}{\omega_{A,k}^{(p)}} \frac{\|A_{I,J} \mathbf{u}\|_p}{\|\mathbf{u}\|_\infty}.$$

In the above derivation, the second inequality follows from $\|(\mathbf{v}/\|\mathbf{v}\|_p)^{p-1}\|_{p/(p-1)} = 1$ and the Hölder's inequality.

The fifth inequality follows from $\|\mathbf{v}\|_p \leq k^{\max(0,1/p-0.5)}\|\mathbf{v}\|_2$ for all $\mathbf{v} \in R^k$, so that

$$\frac{\|A_{I,I}^{-1} A_{I,J}\mathbf{u}\|_p}{\|\mathbf{u}\|_\infty} \leq \frac{\|A_{I,I}^{-1} A_{I,J}\mathbf{u}\|_2}{\|\mathbf{u}\|_\infty} k^{\max(0,1/p-0.5)}.$$

The sixth inequality follows from $\|A_{I,I}^{-1}\mathbf{v}\|_p \leq \|\mathbf{v}\|_p/\mu_{A,k}^{(p)}$ for all $\mathbf{v} \in R^k$, so that

$$\frac{\|A_{I,I}^{-1} A_{I,J}\mathbf{u}\|_p}{\|\mathbf{u}\|_\infty} \leq \frac{1}{\mu_{A,k}^{(p)}} \frac{\|A_{I,J}\mathbf{u}\|_p}{\|\mathbf{u}\|_\infty}.$$

The last inequality is due to

$$\begin{aligned}
\frac{\|A_{I,I}\mathbf{v}\|_p}{\|\mathbf{v}\|_p} &= \frac{\|(\mathbf{v}/\|\mathbf{v}\|_p)^{p-1}\|_{p/(p-1)}\|A_{I,I}\mathbf{v}\|_p}{\|\mathbf{v}\|_p} \\
&\geq \frac{(\mathbf{v}^{p-1})^T A_{I,I}\mathbf{v}}{\|\mathbf{v}\|_p^p} = \frac{(\mathbf{v}^{p-1})^T \operatorname{diag}(A_{I,I})\mathbf{v}}{\|\mathbf{v}\|_p^p} \\
&\quad + \frac{(\mathbf{v}^{p-1})^T (A_{I,I} - \operatorname{diag}(A_{I,I}))\mathbf{v}}{\|\mathbf{v}\|_p^p} \\
&\geq \min_i A_{i,i} - \|(\mathbf{v}/\|\mathbf{v}\|_p)^{p-1}\|_{p/(p-1)} \frac{\|(A_{I,I} - \operatorname{diag}(A_{I,I}))\mathbf{v}\|_p}{\|\mathbf{v}\|_p} \\
&\geq \min_i A_{i,i} - \|A_{I,I} - \operatorname{diag}(A_{I,I})\|_p.
\end{aligned}$$

In the above derivation, Hölder's inequality is used to obtain the first two inequalities. The first equality and the last inequality use the fact that $\|(\mathbf{v}/\|\mathbf{v}\|_p)^{p-1}\|_{p/(p-1)} = 1$.

10.2. *Proof of Proposition 3.2.* Let $B \in R^{d \times d}$ be the off-diagonal part of $A$; that is, $A - B$ is the identity matrix. We have $\sup_{i,j}|B_{i,j}| \leq M_A$. Given any $\mathbf{v} \in R^k$, we have

$$\|B_{I,I}\mathbf{v}\|_p \leq M_A k^{1/p}\|\mathbf{v}\|_1 \leq M_A k\|\mathbf{v}\|_p.$$

This implies that $\|B_{I,I}\|_p \leq M_A k$. Therefore, we have

$$\|A_{I,I}\mathbf{v}\|_p \leq \|\mathbf{v}\|_p(1 + Mk).$$

This proves the first claim. Moreover,

$$(1 - M_A k) \leq 1 - \|B_{I,I}\|_p = \min_i A_{i,i} - \|A_{I,I} - \operatorname{diag}(A_{I,I})\|_p.$$

We thus obtain the second claim from Proposition 3.1.

Now, given $\mathbf{u} \in R^\ell$, since $I \cap J = \varnothing$, we have

$$\|A_{I,J}\mathbf{u}\|_p = \|B_{I,J}\mathbf{u}\|_p \leq M_A k^{1/p}\|\mathbf{u}\|_1 \leq M_A k^{1/p}\ell\|\mathbf{u}\|_\infty.$$

This implies the third claim. The last two claims follow from Proposition 3.1.

10.3. *Some auxiliary results.*

LEMMA 10.1. *Consider $k, \ell > 0$ and $p \in [1, \infty]$. Given any $\beta, \mathbf{v} \in R^d$, let $\beta = \beta_F + \beta_G$, where $\operatorname{supp}_0(\beta_F) \cap \operatorname{supp}_0(\beta_G) = \varnothing$ and $|\operatorname{supp}_0(\beta_F)| = k$. Let $J$ be the indices of the $\ell$ largest components of $\beta_G$ (in absolute values), and $I = \operatorname{supp}_0(\beta_F) \cup J$. If $\operatorname{supp}_0(\mathbf{v}) \subset I$, then*

$$\mathbf{v}^T A \beta \geq \mathbf{v}_I^T A_{I,I} \beta_I - \|A_{I,I} \mathbf{v}_I\|_{p/(p-1)} \gamma_{A,k+\ell,\ell}^{(p)} \|\beta_G\|_1 \ell^{-1},$$

$$\max(0, ((\mathbf{v}/\|\mathbf{v}\|_p)^{p-1})^T A \beta) \geq \omega_{A,k+\ell}^{(p)} (\|\mathbf{v}_I\|_p - \pi_{A,k+\ell,\ell}^{(p)} \ell^{-1} \|\beta_G\|_1)$$
$$- \rho_{A,k+\ell}^{(p)} \|\beta_I - \mathbf{v}_I\|_p.$$

PROOF. In order to prove the first inequality, we may assume, without loss of generality, that $\beta = [\beta_1, \ldots, \beta_d]$ where $\operatorname{supp}_0(\beta_F) = \{1, 2, \ldots, k\}$, and when $j > k$, $\beta_j$ is arranged in descending order of $|\beta_j|$: $|\beta_{k+1}| \geq |\beta_{k+2}| \geq \cdots \geq |\beta_d|$. Let $J_0 = \{1, \ldots, k\}$, and let $J_s = \{k + (s-1)\ell + 1, \ldots, k + s\ell\}$ $(s = 1, 2, \ldots)$, except the largest index in the last block stops at $d$. Note that, in this definition, we require that $J_1 = J$ and $I = J_0 \cup J_1$. We have $\|\beta_{J_s}\|_\infty \leq \|\beta_{J_{s-1}}\|_1 \ell^{-1}$ when $s > 1$, which implies that $\sum_{s>1} \|\beta_{J_s}\|_\infty \leq \|\beta_G\|_1 \ell^{-1}$. This gives

$$\mathbf{v}^T A \beta = \mathbf{v}_I^T A_{I,I} \beta_I + \sum_{s>1} \mathbf{v}_I^T A_{I,J_s} \beta_{J_s}$$

$$\geq \mathbf{v}_I^T A_{I,I} \beta_I - \|A_{I,I} \mathbf{v}_I\|_{p/(p-1)} \sum_{s>1} \|A_{I,I}^{-1} A_{I,J_s} \beta_{J_s}\|_p$$

$$\geq \mathbf{v}_I^T A_{I,I} \beta_I - \gamma_{A,|I|,\ell}^{(p)} \|A_{I,I} \mathbf{v}_I\|_{p/(p-1)} \sum_{s>1} \|\beta_{J_s}\|_\infty$$

$$\geq \mathbf{v}_I^T A_{I,I} \beta_I - \gamma_{A,|I|,\ell}^{(p)} \|A_{I,I} \mathbf{v}_I\|_{p/(p-1)} \|\beta_G\|_1 \ell^{-1}.$$

The first inequality in the above derivation is due to Hölder's inequality. This proves the first inequality of the lemma.

The proof of the second inequality is similar, but with a slightly different estimate. We can assume that the right-hand side is positive (the inequality is trivial otherwise). It implies that $(\mathbf{v}_I^{p-1})^T A_{I,I} \mathbf{v}_I > 0$, since, otherwise, $\omega_{A,|I|}^{(p)} = 0$:

$$(\mathbf{v}^{p-1})^T A \beta = (\mathbf{v}_I^{(p-1)})^T A_{I,I} (\beta_I - \mathbf{v}_I) + (\mathbf{v}_I^{(p-1)})^T A_{I,I} \mathbf{v}_I + \sum_{s>1} (\mathbf{v}_I^{p-1})^T A_{I,J_s} \beta_{J_s}$$

$$\geq (\mathbf{v}_I^{(p-1)})^T A_{I,I} (\beta_I - \mathbf{v}_I)$$

$$+ (\mathbf{v}_I^{(p-1)})^T A_{I,I} \mathbf{v}_I \left[ 1 - \pi_{A,|I|,\ell}^{(p)} \sum_{s>1} \|\beta_{J_s}\|_\infty / \|\mathbf{v}_I\|_p \right]$$

$$\geq -\rho_{A,|I|}^{(p)}\|\mathbf{v}_I^{(p-1)}\|_{p/(p-1)}\|\beta_I - \mathbf{v}_I\|_p$$

$$+ (\mathbf{v}_I^{(p-1)})^T A_{I,I}\mathbf{v}_I[1 - \pi_{A,|I|,\ell}^{(p)}\ell^{-1}\|\beta_G\|_1/\|\mathbf{v}_I\|_p]$$

$$\geq -\rho_{A,|I|}^{(p)}\|\mathbf{v}_I^{(p-1)}\|_{p/(p-1)}\|\beta_I - \mathbf{v}_I\|_p$$

$$+ \omega_{A,|I|}^{(p)}\|\mathbf{v}_I\|_p^p[1 - \pi_{A,|I|,\ell}^{(p)}\ell^{-1}\|\beta_G\|_1/\|\mathbf{v}_I\|_p].$$

The second inequality in the above derivation is due to Hölder's inequality. The last inequality assumes that the right-hand side is nonnegative. Observe that $\|\mathbf{v}_I^{(p-1)}\|_{p/(p-1)} = \|\mathbf{v}_I\|_p^{p-1}$; thus, we obtain the second inequality of the lemma.  $\square$

Lemma 10.2.   *Consider the decomposition of any target vector $\bar{\beta} = \bar{\beta}_F + \bar{\beta}_G$ such that $\{1, 2, \ldots, d\} = F \cup G$ and $F \cap G = \varnothing$. Consider the solution $\hat{\beta}$ to the following, more general, problem instead of (1):*

$$(5) \qquad \hat{\beta} = \underset{\beta \in R^d}{\arg\min}\left[\frac{1}{n}\sum_{i=1}^n (\beta^T \mathbf{x}_i - \mathbf{y}_i)^2 + \lambda\sum_{j \notin \hat{F}}|\beta_j|\right],$$

*where $\hat{F} \subset F$. Let $\Delta\hat{\beta} = \hat{\beta} - \bar{\beta}$, $\hat{A} = \frac{1}{n}\sum_{i=1}^n \mathbf{x}_i\mathbf{x}_i^T$ and $\hat{\varepsilon} = \frac{1}{n}\sum_{i=1}^n (\bar{\beta}^T\mathbf{x}_i - \mathbf{y}_i)\mathbf{x}_i$. If we pick a sufficiently large $\lambda$ in (5) such that $\lambda > 2\|\hat{\varepsilon}\|_\infty$, then*

$$\|\hat{\beta}_G\|_1 \leq \frac{2\|\hat{\varepsilon}\|_\infty + \lambda}{\lambda - 2\|\hat{\varepsilon}\|_\infty}(\|\Delta\hat{\beta}_F\|_1 + \|\bar{\beta}_G\|_1).$$

Proof.   We define the derivative of $\|\beta\|_1$ as $\mathrm{sgn}(\beta)$, where, for $\beta = [\beta_1, \ldots, \beta_d] \in R^d$, $\mathrm{sgn}(\beta) = [\mathrm{sgn}(\beta_1), \ldots, \mathrm{sgn}(\beta_d)] \in R^d$ is defined as $\mathrm{sgn}(\beta_j) = 1$ when $\beta_j > 0$, $\mathrm{sgn}(\beta_j) = -1$ when $\beta_j < 0$ and $\mathrm{sgn}(\beta_j) \in [-1, 1]$ when $\beta_j = 0$. We start with the first-order condition

$$\frac{2}{n}\sum_{i=1}^n (\hat{\beta}^T\mathbf{x}_i - \mathbf{y}_i)\mathbf{x}_i + \lambda g(\hat{\beta}) = 0,$$

where $g(\hat{\beta}) = [g(\hat{\beta}_1), \ldots, g(\hat{\beta}_d)]$, with $g(\hat{\beta}_j) = 0$ when $j \in \hat{F}$ and $g(\hat{\beta}_j) = \mathrm{sgn}(\hat{\beta}_j)$, otherwise. This implies that

$$2\hat{A}\Delta\hat{\beta} + \lambda g(\hat{\beta}) = -\frac{2}{n}\sum_{i=1}^n (\bar{\beta}^T\mathbf{x}_i - \mathbf{y}_i)\mathbf{x}_i.$$

Therefore, for all $\mathbf{v} \in R^d$, we have

$$(6) \qquad 2\mathbf{v}^T\hat{A}\Delta\hat{\beta} \leq -2\mathbf{v}^T\hat{\varepsilon} - \lambda\mathbf{v}^T g(\hat{\beta}).$$

Now, let $\mathbf{v} = \Delta\hat{\beta}$ in (6), and use the fact that $\hat{\beta}_G^T g(\hat{\beta}_G) = \hat{\beta}_G^T \operatorname{sgn}(\hat{\beta}_G) = \|\hat{\beta}_G\|_1$ as well as $\|g(\hat{\beta})\|_\infty \leq 1$. We obtain

$$0 \leq 2\Delta\hat{\beta}^T \hat{A}\Delta\hat{\beta} \leq 2|\Delta\hat{\beta}^T \hat{\varepsilon}| - \lambda\Delta\hat{\beta}^T g(\hat{\beta})$$

$$\leq 2\|\Delta\hat{\beta}\|_1\|\hat{\varepsilon}\|_\infty - \lambda\Delta\hat{\beta}_F^{\ T} g(\hat{\beta}) - \lambda\hat{\beta}_G^T g(\hat{\beta}) + \lambda\bar{\beta}_G^T g(\hat{\beta})$$

$$\leq 2(\|\Delta\hat{\beta}_F\|_1 + \|\hat{\beta}_G\|_1 + \|\bar{\beta}_G\|_1)\|\hat{\varepsilon}\|_\infty + \lambda\|\Delta\hat{\beta}_F\|_1 - \lambda\|\hat{\beta}_G\|_1 + \lambda\|\bar{\beta}_G\|_1$$

$$= (2\|\hat{\varepsilon}\|_\infty - \lambda)\|\hat{\beta}_G\|_1 + (2\|\hat{\varepsilon}\|_\infty + \lambda)(\|\Delta\hat{\beta}_F\|_1 + \|\bar{\beta}_G\|_1).$$

By rearranging the above inequality, we obtain the desired bound.  $\square$

LEMMA 10.3.  *Let the conditions of Lemma 10.2 hold. Let $J$ be the indices of the largest $\ell$ coefficients (in absolute value) of $\Delta\hat{\beta}_G$ and $I = F \cup J$. If $\lambda \geq 4(2-t)t^{-1}\|\hat{\varepsilon}\|_\infty$ for some $t \in (0,1)$, then $\forall p \in [1, \infty]$,*

$$\|\Delta\hat{\beta}\|_1 \leq 4k^{1-1/p}\|\Delta\hat{\beta}_I\|_p + 4\|\bar{\beta}_G\|_1,$$

$$\|\Delta\hat{\beta}\|_p \leq (1 + 3(k/\ell)^{1-1/p})\|\Delta\hat{\beta}_I\|_p + 4\|\bar{\beta}_G\|_1\ell^{1/p-1}.$$

PROOF.  The condition on $\lambda$ implies that $(\lambda + 2\|\hat{\varepsilon}\|_\infty)/(\lambda - 2\|\hat{\varepsilon}\|_\infty) \leq (4-t)/(4-3t) \leq 3$. We have, from Lemma 10.2,

$$\|\Delta\hat{\beta}_G\|_1 \leq \|\bar{\beta}_G\|_1 + \|\hat{\beta}_G\|_1 \leq 3\|\Delta\hat{\beta}_F\|_1 + 4\|\bar{\beta}_G\|_1.$$

Therefore, $\|\Delta\hat{\beta} - \Delta\hat{\beta}_I\|_\infty \leq \|\Delta\hat{\beta}_G\|_1/\ell \leq (3\|\Delta\hat{\beta}_F\|_1 + 4\|\bar{\beta}_G\|_1)/\ell$, which implies that

$$\|\Delta\hat{\beta} - \Delta\hat{\beta}_I\|_p \leq (\|\Delta\hat{\beta}_G\|_1\|\Delta\hat{\beta} - \Delta\hat{\beta}_I\|_\infty^{p-1})^{1/p} \leq (3\|\Delta\hat{\beta}_F\|_1 + 4\|\bar{\beta}_G\|_1)\ell^{1/p-1}.$$

Now, the first inequality in the proof also implies that

$$\|\Delta\hat{\beta}\|_1 \leq 4\|\Delta\hat{\beta}_F\|_1 + 4\|\bar{\beta}_G\|_1.$$

By combining the previous two inequalities with $\|\Delta\hat{\beta}_F\|_1 \leq k^{1-1/p}\|\Delta\hat{\beta}_I\|_p$, we obtain the desired bounds.  $\square$

LEMMA 10.4.  *Let the conditions of Lemma 10.2 hold. Let $J$ be the indices of the largest $\ell$ coefficients (in absolute value) of $\Delta\hat{\beta}_G$ and $I = F \cup J$. Assume that $\hat{A}_{I,I}$ is invertible. If $t = 1 - \pi_{\hat{A},k+\ell,\ell}^{(p)}k^{1-1/p}\ell^{-1} > 0$ and $\lambda \geq 4(2-t)t^{-1}\|\hat{\varepsilon}\|_\infty$, then*

$$\|\Delta\hat{\beta}_I\|_p \leq \frac{2}{t\omega_{\hat{A},k+\ell}^{(p)}}[\rho_{\hat{A},k+\ell}^{(p)}\|\bar{\beta}_G\|_p + (k - |F_0|)^{1/p}\lambda + \|\hat{\varepsilon}_{F_0}\|_p]$$

$$+ \frac{8\pi_{\hat{A},k+\ell,\ell}^{(p)}\|\bar{\beta}_G\|_1/\ell}{t} + \|\bar{\beta}_G\|_p,$$

*where $F_0$ is any subset of $\hat{F}$.*

Proof.  The condition of $\lambda$ implies that $(\lambda + 2\|\hat{\varepsilon}\|_\infty)/(\lambda - 2\|\hat{\varepsilon}\|_\infty) \leq (4-t)/(4-3t)$. Therefore, if we let $\Delta\dot{\beta} = \Delta\dot{\beta}_I = \Delta\dot{\beta}_I + \bar{\beta}_J$, then

$$\max(0, ((\Delta\dot{\beta}_I/\|\Delta\dot{\beta}_I\|_p)^{p-1})^T \hat{A}\Delta\hat{\beta}) + \rho^{(p)}_{\hat{A},k+\ell}\|\bar{\beta}_J\|_p$$

$$\geq \omega^{(p)}_{\hat{A},k+\ell}(\|\Delta\dot{\beta}_I\|_p - \pi^{(p)}_{\hat{A},k+\ell,\ell}\ell^{-1}\|\Delta\hat{\beta}_G\|_1)$$

$$\geq \omega^{(p)}_{\hat{A},k+\ell}[\|\Delta\dot{\beta}_I\|_p$$
$$- \pi^{(p)}_{\hat{A},k+\ell,\ell}\ell^{-1}((4-t)(4-3t)^{-1}(\|\Delta\hat{\beta}_F\|_1 + \|\bar{\beta}_G\|_1) + \|\bar{\beta}_G\|_1)]$$

$$\geq \omega^{(p)}_{\hat{A},k+\ell}(\|\Delta\dot{\beta}_I\|_p - (1-t)(4-t)(4-3t)^{-1}\|\Delta\dot{\beta}_I\|_p$$
$$- 4\pi^{(p)}_{\hat{A},k+\ell,\ell}\ell^{-1}\|\bar{\beta}_G\|_1)$$

$$\geq 0.5t\omega^{(p)}_{\hat{A},k+\ell}\|\Delta\dot{\beta}_I\|_p - 4\omega^{(p)}_{\hat{A},k+\ell}\pi^{(p)}_{\hat{A},k+\ell,\ell}\ell^{-1}\|\bar{\beta}_G\|_1.$$

In the above derivation, the first inequality is due to Lemma 10.1, and the second inequality is due to Lemma 10.2. The third inequality uses $\|\Delta\hat{\beta}_F\|_1 = \|\Delta\hat{\beta}_F\|_1 \leq k^{1-1/p}\|\Delta\dot{\beta}_I\|_p$ and $(4-t)/(4-3t) \leq 3$. The last inequality follows from $1 - (1-t)(4-t)(4-3t)^{-1} \geq 0.5t$.

If $(\Delta\dot{\beta}_I^{p-1})^T\hat{A}\Delta\hat{\beta} \leq 0$, then the above inequality, together with $\|\Delta\hat{\beta}_I\|_p \leq \|\Delta\dot{\beta}_I\|_p + \|\bar{\beta}_J\|_p \leq \|\Delta\dot{\beta}_I\|_p + \|\bar{\beta}_G\|_p$, already implies the lemma. Therefore in the following, we can assume that

$$((\Delta\dot{\beta}_I/\|\Delta\dot{\beta}_I\|_p)^{p-1})^T\hat{A}\Delta\hat{\beta} \geq 0.5t\omega^{(p)}_{\hat{A},k+\ell}\|\Delta\dot{\beta}_I\|_p$$
$$- 4\omega^{(p)}_{\hat{A},k+\ell}\pi^{(p)}_{\hat{A},k+\ell,\ell}\ell^{-1}\|\bar{\beta}_G\|_1 - \rho^{(p)}_{\hat{A},k+\ell}\|\bar{\beta}_J\|_p.$$

Moreover, we obtain, from (6) with $\mathbf{v} = \Delta\dot{\beta}_I^{p-1}$, the following:

$$(\Delta\dot{\beta}_I^{p-1})^T\hat{A}\Delta\hat{\beta}$$

$$\leq |(\Delta\dot{\beta}_I^{p-1})^T\hat{\varepsilon}| - \lambda(\Delta\dot{\beta}_I^{p-1})^T g(\hat{\beta})/2$$

$$\leq |(\hat{\beta}_J^{p-1})^T\hat{\varepsilon}| - \lambda(\hat{\beta}_J^{p-1})^T g(\hat{\beta})/2 + |(\Delta\hat{\beta}_{F-F_0}^{p-1})^T\hat{\varepsilon}|$$
$$- \lambda(\Delta\hat{\beta}_{F-F_0}^{p-1})^T g(\hat{\beta})/2 + |(\Delta\hat{\beta}_{F_0}^{p-1})^T\hat{\varepsilon}|$$

$$\leq (\|\hat{\varepsilon}\|_\infty - \lambda/2)\|\hat{\beta}_J^{p-1}\|_1 + (\|\hat{\varepsilon}\|_\infty + \lambda/2)\|\Delta\hat{\beta}_{F-F_0}^{p-1}\|_1 + |(\Delta\hat{\beta}_{F_0}^{p-1})^T\hat{\varepsilon}|$$

$$\leq \lambda\|\Delta\hat{\beta}_{F-F_0}^{p-1}\|_1 + |(\Delta\hat{\beta}_{F_0}^{p-1})^T\hat{\varepsilon}|$$

$$\leq (k-|F_0|)^{1/p}\lambda\|\Delta\hat{\beta}_{F-F_0}^{p-1}\|_{p/(p-1)} + \|\Delta\hat{\beta}_{F_0}^{p-1}\|_{p/(p-1)}\|\hat{\varepsilon}_{F_0}\|_p$$

$$\leq ((k-|F_0|)^{1/p}\lambda + \|\hat{\varepsilon}_{F_0}\|_p)\|\Delta\dot{\beta}_I\|_p^{p-1}.$$

In the above derivation, the second inequality uses $g(\hat{\beta}_{F_0}) = 0$; the third inequality uses the fact that $\|g(\hat{\beta})\|_\infty \leq 1$ and $(\hat{\beta}_J^{p-1})^T g(\hat{\beta}) = \|\hat{\beta}_J^{p-1}\|_1$; the fourth inequality uses $\|\hat{\varepsilon}\|_\infty \leq 0.5\lambda$; and the last inequality uses the fact that $\|\beta^{p-1}\|_{p/(p-1)} = \|\beta\|_p^{p-1}$. Now, by combining the above two estimates, together with $\|\Delta\hat{\beta}_I\|_p \leq \|\Delta\dot{\beta}_I\|_p + \|\bar{\beta}_J\|_p \leq \|\Delta\dot{\beta}_I\|_p + \|\bar{\beta}_G\|_p$, we obtain the desired bound. $\square$

LEMMA 10.5. *Let the conditions of Lemma 10.2 hold. Let $J$ be the indices of the largest $\ell$ coefficients (in absolute value) of $\Delta\hat{\beta}_G$, and $I = F \cup J$. Assume that $\hat{A}_{I,I}$ is invertible. Let $p \in [1, \infty]$. If $t = 1 - \gamma_{\hat{A}, k+\ell, \ell}^{(p)} k^{1-1/p} \ell^{-1} > 0$ and $\lambda \geq 4(2-t)t^{-1}\|\hat{\varepsilon}\|_\infty$. Then,*

$$\|\Delta\hat{\beta}_I\|_p \leq \frac{2}{t}[4\gamma_{\hat{A}, k+\ell, \ell}^{(p)} \ell^{-1} \|\bar{\beta}_G\|_1 + \lambda(k+\ell)^{1/p}/\mu_{\hat{A}, k+\ell}^{(p)}].$$

PROOF. Consider $\mathbf{v} \in R^d$ such that $\text{supp}_0(\mathbf{v}) \subset I$. We have, from Lemma 10.1 and (6),

$$\mathbf{v}^T \hat{A}\Delta\hat{\beta}_I - \|\hat{A}_{I,I}\mathbf{v}_I\|_{p/(p-1)}\gamma_{\hat{A}, k+\ell, \ell}^{(p)} \|\Delta\hat{\beta}_G\|_1 \ell^{-1} \leq \mathbf{v}^T \hat{A}\Delta\hat{\beta}$$
$$\leq -\mathbf{v}^T(\hat{\varepsilon} + 0.5\lambda g(\hat{\beta})).$$

Take $\mathbf{v}$ such that $\|\hat{A}_{I,I}\mathbf{v}_I\|_{p/(p-1)} = 1$ and $\mathbf{v}^T \hat{A}\Delta\hat{\beta}_I = \|\Delta\hat{\beta}_I\|_p$. We obtain

$$\|\Delta\hat{\beta}_I\|_p - \gamma_{\hat{A}, k+\ell, \ell}^{(p)}(\|\hat{\beta}_G\|_1 + \|\bar{\beta}_G\|_1)\ell^{-1} \leq \|\hat{A}_{I,I}^{-1}(\hat{\varepsilon}_I + 0.5\lambda g(\hat{\beta}_I))\|_p.$$

By using Lemma 10.2, $(\lambda + 2\|\hat{\varepsilon}\|_\infty)/(\lambda - 2\|\hat{\varepsilon}\|_\infty) \leq (4-t)/(4-3t) \leq 3$ and $1 - (1-t)(4-t)/(4-3t) \geq 0.5t$; thus, we obtain

$$\|\Delta\hat{\beta}_I\|_p - \gamma_{\hat{A}, k+\ell, \ell}^{(p)} \ell^{-1}\|\hat{\beta}_G\|_1$$
$$\geq \|\Delta\hat{\beta}_I\|_p - \gamma_{\hat{A}, k+\ell, \ell}^{(p)} \ell^{-1}(4-t)(4-3t)^{-1}(\|\Delta\hat{\beta}_F\|_1 + \|\bar{\beta}_G\|_1)$$
$$\geq \|\Delta\hat{\beta}_I\|_p - \gamma_{\hat{A}, k+\ell, \ell}^{(p)} \ell^{-1}(4-t)(4-3t)^{-1}(k^{1-1/p}\|\Delta\hat{\beta}_F\|_p + \|\bar{\beta}_G\|_1)$$
$$\geq \|\Delta\hat{\beta}_I\|_p - (1-t)(4-t)(4-3t)^{-1}\|\Delta\hat{\beta}_I\|_p - 3\gamma_{\hat{A}, k+\ell, \ell}^{(p)} \ell^{-1}\|\bar{\beta}_G\|_1$$
$$\geq 0.5t\|\Delta\hat{\beta}_I\|_p - 3\gamma_{\hat{A}, k+\ell, \ell}^{(p)} \ell^{-1}\|\bar{\beta}_G\|_1.$$

Combine the previous two inequalities. We obtain

$$0.5t\|\Delta\hat{\beta}_I\|_p \leq 4\gamma_{\hat{A}, k+\ell, \ell}^{(p)}\|\bar{\beta}_G\|_1\ell^{-1} + \|\hat{A}_{I,I}^{-1}(\hat{\varepsilon}_I + 0.5\lambda g(\hat{\beta}_I))\|_p$$
$$\leq 4\gamma_{\hat{A}, k+\ell, \ell}^{(p)}\|\bar{\beta}_G\|_1\ell^{-1} + (k+\ell)^{1/p}\|\hat{\varepsilon}_I + 0.5\lambda g(\hat{\beta}_I)\|_\infty/\mu_{\hat{A}, k+\ell}^{(p)}.$$

Since $\|\hat{\varepsilon}_I + 0.5\lambda g(\hat{\beta}_I)\|_\infty \leq \lambda$, we obtain the desired bound. $\square$

PROPOSITION 10.1. *Consider $n$ independent random variables $\xi_1, \ldots, \xi_n$ such that $\mathbf{E}e^{t(\xi_i - \mathbf{E}\xi_i)} \leq e^{\sigma_i^2 t^2/2}$ for all $t$ and $i$, then $\forall \varepsilon > 0$:*

$$P\left(\left|\sum_{i=1}^n \xi_i - \sum_{i=1}^n \mathbf{E}\xi_i\right| \geq n\varepsilon\right) \leq 2e^{-n^2\varepsilon^2/(2\sum_{i=1}^n \sigma_i^2)}.$$

PROOF. Let $s_n = \sum_{i=1}^n (\xi_i - \mathbf{E}\xi_i)$; then, by assumption, $\mathbf{E}(e^{ts_n} + e^{-ts_n}) \leq 2e^{\sum_i \sigma_i^2 t^2/2}$, which implies that $P(|s_n| \geq n\varepsilon)e^{tn\varepsilon} \leq 2e^{\sum_i \sigma_i^2 t^2/2}$. Now, let $t = n\varepsilon / \sum_i \sigma_i^2$; thus, we obtain the desired bound. $\square$

PROPOSITION 10.2. *Consider $n$ independent random variables $\xi_1, \ldots, \xi_n$, such that $\mathbf{E}\xi_i = 0$ and $\mathbf{E}e^{t\xi_i} \leq e^{\sigma^2 t^2/2}$ for all $t$ and $i$. Let $\mathbf{z}_1, \ldots, \mathbf{z}_n \in R^d$ be $n$ fixed vectors, and let $a_n = (\sum_{i=1}^n \|\mathbf{z}_i\|_2^2)^{1/2}$. Then, $\forall \varepsilon > 0$:*

$$P\left(\left\|\sum_{i=1}^n \xi_i \mathbf{z}_i\right\|_2 \geq a_n(\sigma + \varepsilon)\right) \leq e^{-\varepsilon^2/(20\sigma^2)}.$$

PROOF. For each $i$, let $\xi_i'$ be an identically distributed and independent copy of $\xi_i$ and $h(\cdot)$ be any real-valued function such that $h(\xi_i) - h(\xi_i') \leq |\xi_i| + |\xi_i'|$. Then,

$$
\begin{aligned}
&\mathbf{E}_{\xi_i} e^{t(h(\xi_i) - \mathbf{E}_{\xi_i'} h(\xi_i'))} \\
&= 1 + \sum_{k=2}^\infty \frac{t^k}{k!} \mathbf{E}_{\xi_i}(h(\xi_i) - \mathbf{E}_{\xi_i'} h(\xi_i'))^k \\
&\leq 1 + \sum_{k=2}^\infty \frac{t^k}{k!} \mathbf{E}_{\xi_i}(|\xi_i| + \mathbf{E}_{\xi_i'}|\xi_i'|)^k \leq 1 + \sum_{k=2}^\infty \frac{(2t)^k}{k!} \mathbf{E}_{\xi_i}|\xi_i|^k \\
&= 1 + \sum_{k=1}^\infty \left[\frac{1}{(2k)!}\mathbf{E}_{\xi_i}|2t\xi_i|^{2k} + \frac{1}{(2k+1)!}\mathbf{E}_{\xi_i}|2t\xi_i|^{2k+1}\right] \\
&\leq 1 + \sum_{k=1}^\infty \left[\frac{1}{(2k)!}\mathbf{E}|2t\xi_i|^{2k} + \frac{0.5}{(2k)!}\mathbf{E}|2t\xi_i|^{2k} + \frac{1}{(2k+2)!}\mathbf{E}|2t\xi_i|^{2k+2}\right] \\
&\leq 1 + 2.5 \sum_{k=1}^\infty \frac{1}{(2k)!}\mathbf{E}|2t\xi_i|^{2k} = 1 + 1.25(\mathbf{E}e^{2t\xi_i} + \mathbf{E}e^{-2t\xi_i} - 2) \\
&\leq 1 + 1.25(2e^{2t^2\sigma^2} - 2) \leq e^{5t^2\sigma^2}.
\end{aligned}
$$

The second inequality is due to Jensen's inequality. In the third inequality, we have used $|a|^{2k+1}/(2k+1)! \leq 0.5|a|^{2k}/(2k)! + |a|^{2k+2}/(2k+2)!$. The last inequality can be obtained by comparing the Taylor expansion of the function $e^x$ on both sides.

Now, let $s_j = \mathbf{E}_{\xi_{j+1},\ldots,\xi_n} \|\sum_{i=1}^n \xi_i \mathbf{z}_i\|_2$. If we regard $h(\xi_j) = s_j/\|\mathbf{z}_j\|_2$ as a function of $\xi_j$ (with variables $\xi_1,\ldots,\xi_{j-1}$ fixed), then $s_j - s_{j-1} = (h(\xi_j) - \mathbf{E}_{\xi'_j} h(\xi'_j))\|\mathbf{z}_j\|_2$ and $h(\xi_j) - h(\xi'_j) \leq |\xi_j| + |\xi'_j|$. Therefore, from the above inequality, we have $\mathbf{E}_{\xi_j} e^{t(s_j - s_{j-1})} \leq e^{5\|\mathbf{z}_j\|_2^2 t^2 \sigma^2}$ and

$$\mathbf{E}_{\xi_1,\ldots,\xi_j} e^{ts_j} = \mathbf{E}_{\xi_1,\ldots,\xi_{j-1}} e^{ts_{j-1}} \mathbf{E}_{\xi_j} e^{t(s_j - s_{j-1})} \leq e^{5\|\mathbf{z}_j\|_2^2 \sigma^2 t^2} \mathbf{E}_{\xi_1,\ldots,\xi_{j-1}} e^{ts_{j-1}}.$$

By induction, we obtain $\mathbf{E}_{\xi_1,\ldots,\xi_n} e^{ts_n} \leq e^{5\sigma^2 t^2 a_n^2} e^{ts_0}$, which implies that $P(s_n \geq s_0 + a_n \varepsilon) e^{t(s_0 + a_n \varepsilon)} \leq e^{5a_n^2 \sigma^2 t^2} e^{ts_0}$. Let $t = \varepsilon/(10 a_n \sigma^2)$, we have $P(s_n \geq s_0 + a_n \varepsilon) \leq e^{-\varepsilon^2/(20\sigma^2)}$.

Note that

$$\mathbf{E}\xi_i^2 = \lim_{t \to 0} \frac{2}{t^2}(\mathbf{E}_{\xi_i} e^{t\xi_i} - 1) \leq \lim_{t \to 0} \frac{2(e^{\sigma^2 t^2/2} - 1)}{t^2} = \sigma^2.$$

Therefore, $s_0 = \mathbf{E}\|\sum_{i=1}^n \xi_i \mathbf{z}_i\|_2 \leq (\sum_{i=1}^n \mathbf{E}\xi_i^2 \|\mathbf{z}_i\|_2^2)^{1/2} \leq a_n \sigma$. This leads to the desired bound. $\square$

10.4. *Proof of Theorem 4.1.* Let $F$ be the indices corresponding to the largest $k$ coefficients of $\bar{\beta}$ in absolute value. We only need to estimate $\|\hat{\varepsilon}\|_\infty$ and then apply Lemmas 10.4, 10.5 and 10.3, with $F_0 = \varnothing$. By Proposition 10.1, we have

$$P\left[\sup_j \sum_{i=1}^n \mathbf{x}_{i,j}^2 \leq na^2 \text{ and } \sup_j \left|\frac{1}{n}\sum_{i=1}^n (\mathbf{y}_i - \mathbf{E}\mathbf{y}_i)\mathbf{x}_{i,j}\right| \geq \varepsilon\right]$$

$$\leq d \sup_j P\left[\left|\frac{1}{n}\sum_{i=1}^n (\mathbf{y}_i - \mathbf{E}\mathbf{y}_i)\mathbf{x}_{i,j}\right| \geq \varepsilon \,\Big|\, \sup_{j'} \sum_{i=1}^n \mathbf{x}_{i,j'}^2 \leq na^2\right]$$

$$\leq 2d \sup_j e^{-n^2 \varepsilon^2/(2\sigma^2 \sum_{i=1}^n \mathbf{x}_{i,j}^2)} \qquad \left(\text{subject to } \sup_{j'} \sum_{i=1}^n \mathbf{x}_{i,j'}^2 \leq na^2\right)$$

$$\leq 2d e^{-n\varepsilon^2/(2\sigma^2 a^2)}.$$

Therefore, with probability larger than $1 - \delta$, if $\sup_j \sum_{i=1}^n \mathbf{x}_{i,j}^2 \leq na^2$, then

$$\sup_j \left|\frac{1}{n}\sum_{i=1}^n (\mathbf{y}_i - \mathbf{E}\mathbf{y}_i)\mathbf{x}_{i,j}\right| \leq \sigma a \sqrt{\frac{2\ln(2d/\delta)}{n}}.$$

The latter implies that

$$\|\hat{\varepsilon}\|_\infty = \sup_j \left|\frac{1}{n}\sum_{i=1}^n (\bar{\beta}^T \mathbf{x}_i - \mathbf{y}_i)\mathbf{x}_{i,j}\right| \leq \sigma a \sqrt{\frac{2\ln(2d/\delta)}{n}} + \left\|\frac{1}{n}\sum_{i=1}^n (\bar{\beta}^T \mathbf{x}_i - \mathbf{E}\mathbf{y}_i)\mathbf{x}_i\right\|_\infty.$$

With this bound, the condition of Lemma 10.3 is satisfied. Using $k/\ell \leq 1$, we obtain, for $q = 1, p$,

$$\|\Delta\hat{\beta}\|_q \leq 4k^{1/q-1/p}\|\Delta\hat{\beta}_I\|_p + 4\|\bar{\beta}_G\|_1 \ell^{1/q-1}.$$

This estimate, together with Lemma 10.4 (let $F_0 = \varnothing$), leads to the first claim of the theorem; with Lemma 10.5, it gives the second claim.

10.5. *Proof of Corollary 4.1.* The condition $M_{\hat{A}}(k+\ell) \leq (1-t)/(2-t)$ is equivalent to $t \leq 1 - M_{\hat{A}}(k+\ell)/(1 - M_{\hat{A}}(k+\ell))$. It implies $t \leq 1 - M_{\hat{A}}(k+\ell)^{1/p}k^{1-1/p}/(1 - M_{\hat{A}}(k+\ell))$. Now, using Proposition 3.2, it implies that the condition $t \leq 1 - \pi^{(p)}_{\hat{A},k+\ell,\ell}k^{1-1/p}/\ell$ in the first claim of Theorem 4.1 is satisfied. Therefore, we have (with $q = p$)

$$\|\Delta\hat{\beta}\|_p \leq \frac{8}{t\omega^{(p)}_{\hat{A},k+\ell}}[\rho^{(p)}_{\hat{A},k+\ell}r^{(p)}_k(\bar{\beta}) + k^{1/p}\lambda] + 4r^{(p)}_k(\bar{\beta})$$

$$+ \left[\frac{32}{t}\pi^{(p)}_{\hat{A},k+\ell,\ell}\ell^{-1} + 4\ell^{1/p-1}\right]r^{(1)}_k(\bar{\beta}).$$

Now, using Proposition 3.2 again, the inequality can be simplified to

$$\|\Delta\hat{\beta}\|_p \leq \frac{8}{t(1 - M_{\hat{A}}(k+\ell))}[(1 + M_{\hat{A}}(k+\ell))r^{(p)}_k(\bar{\beta}) + k^{1/p}\lambda]$$

$$+ \frac{32M_{\hat{A}}(k+\ell)^{1/p}}{t(1 - M_{\hat{A}}(k+\ell))}r^{(1)}_k(\bar{\beta}) + 4r^{(p)}_k(\bar{\beta}) + 4r^{(1)}_k(\bar{\beta})\ell^{1/p-1}.$$

Now, by using the condition $M_{\hat{A}}(k+\ell) \leq (1-t)/(2-t) \leq 0.5$ to eliminate $M_{\hat{A}}$ and simplify the result, we obtain the desired bound.

10.6. *Proof of Proposition 5.1.* We construct a sequence $\beta^{(k)}$ with a greedy algorithm as follows. Let $\beta^{(0)} = \bar{\beta}$, and, for $k = 1, 2, \ldots$, we perform the following steps:

- $j^{(k)} = \arg\max_j |\sum_{i=1}^n (\beta^{(k-1)T}\mathbf{x}_i - \mathbf{E}\mathbf{y}_i)\mathbf{x}_{i,j}|$;
- $\alpha^{(k)} = -\sum_{i=1}^n (\beta^{(k-1)T}\mathbf{x}_i - \mathbf{E}\mathbf{y}_i)\mathbf{x}_{i,j^{(k)}}/(na^2)$;
- $\beta^{(k)} = \beta^{(k-1)} + \alpha^{(k)}\mathbf{e}_{j^{(k)}}$, where $\mathbf{e}_j \in R^d$ is the vector of zeros, except for the $j$th component being one.

The following derivation holds for the above procedure:

$$\sum_{i=1}^n (\beta^{(k)T}\mathbf{x}_i - \mathbf{E}\mathbf{y}_i)^2$$

$$= \sum_{i=1}^n (\beta^{(k-1)T}\mathbf{x}_i - \mathbf{E}\mathbf{y}_i)^2 + \alpha^{(k)2}\sum_{i=1}^n \mathbf{x}^2_{i,j^{(k)}}$$

$$+ 2\alpha^{(k)} \sum_{i=1}^{n} (\beta^{(k-1)T}\mathbf{x}_i - \mathbf{E}\mathbf{y}_i)\mathbf{x}_{i,j^{(k)}}$$

$$\leq \sum_{i=1}^{n} (\beta^{(k-1)T}\mathbf{x}_i - \mathbf{E}\mathbf{y}_i)^2 - \left| \sum_{i=1}^{n} (\beta^{(k-1)T}\mathbf{x}_i - \mathbf{E}\mathbf{y}_i)\mathbf{x}_{i,j^{(k)}} \right|^2 \Big/ (na^2)$$

$$= \sum_{i=1}^{n} (\beta^{(k-1)T}\mathbf{x}_i - \mathbf{E}\mathbf{y}_i)^2 - \left\| \sum_{i=1}^{n} (\beta^{(k-1)T}\mathbf{x}_i - \mathbf{E}\mathbf{y}_i)\mathbf{x}_i \right\|_\infty^2 \Big/ (na^2).$$

In the above derivation, first equality is simple algebra; the inequality uses the definition of $a$ and $\alpha^{(k)}$; and the last equality uses the definition of $j^{(k)}$. Since $\sum_{i=1}^{n} (\beta^{(k+1)T}\mathbf{x}_i - \mathbf{E}\mathbf{y}_i)^2 \geq 0$, we obtain

$$\sum_{k'=0}^{k} \left\| \sum_{i=1}^{n} (\beta^{(k')T}\mathbf{x}_i - \mathbf{E}\mathbf{y}_i)\mathbf{x}_i \right\|_\infty^2 \leq na^2 \sum_{i=1}^{n} (\beta^{(0)T}\mathbf{x}_i - \mathbf{E}\mathbf{y}_i)^2.$$

Therefore, there exists $k' \leq k$ such that the displayed equation of the proposition holds with $\bar{\beta}^{(k)} = \beta^{(k')}$. Moreover,

$$\sqrt{\mu_{\hat{A},k}^{(2)}} \|(\bar{\beta}^{(k)} - \bar{\beta})\|_2 \leq \|\hat{A}^{1/2}(\bar{\beta}^{(k)} - \bar{\beta})\|_2$$

$$\leq \|\hat{A}^{1/2}(\bar{\beta}^{(k)} - \mathbf{E}\mathbf{y})\|_2 + \|\hat{A}^{1/2}(\bar{\beta} - \mathbf{E}\mathbf{y})\|_2$$

$$\leq 2\|\hat{A}^{1/2}(\bar{\beta} - \mathbf{E}\mathbf{y})\|_2.$$

The proof is complete.

10.7. *Proof of Corollary 5.1.* The proof is just a straightforward application of Corollary 4.1, in which we replace $\bar{\beta}$ by $\bar{\beta}^{(k)}$ of Proposition 5.1 and then replace $k$ by $2k$. This leads to the bound

$$\|\hat{\beta} - \bar{\beta}\|_2 \leq \|\hat{\beta} - \bar{\beta}^{(k)}\|_2 + \|\bar{\beta} - \bar{\beta}^{(k)}\|_2$$

$$\leq \frac{8(2-t)}{t}[1.5r_{2k}^{(2)}(\bar{\beta}^{(k)}) + (2k)^{1/2}\lambda] + 4r_{2k}^{(2)}(\bar{\beta}^{(k)})$$

$$+ \frac{4(8-7t)}{t}r_{2k}^{(1)}(\bar{\beta}^{(k)})\ell^{-1/2} + 2\varepsilon/\sqrt{\mu_{\hat{A},k}^{(2)}}.$$

Note that, from Proposition 3.2, we have $1/\sqrt{\mu_{\hat{A},k}^{(2)}} \leq (1 - kM_{\hat{A}})^{-1/2} < 2$. Moreover, since $\mathrm{supp}_0(\bar{\beta}^{(k)} - \bar{\beta}) \leq k$, we have $r_{2k}^{(p)}(\bar{\beta}^{(k)}) \leq r_k^{(p)}(\bar{\beta})$. This leads to the desired bound.

10.8. *Proof of Corollary 6.1.* Note that Proposition 3.1 implies that $\pi_{\hat{A},k+\ell,\ell}^{(2)} \le \theta_{\hat{A},k+\ell,\ell}^{(2)}/\omega_{\hat{A},k+\ell}^{(2)} \le \rho_{\hat{A},k+\ell}^{(2)}\sqrt{\ell}/\omega_{\hat{A},k+\ell}^{(2)}$. Also note that $\omega_{A,k+\ell}^{(2)} = \mu_{A,k+\ell}^{(2)}$. The first statement of Theorem 4.1 (with $q = p = 2$) implies the desired bound.

10.9. *Proof of Theorem 7.1.* Under the conditions of this theorem, we obtain, from Theorem 4.1, that, with probability $1 - \delta$, the following two claims hold (with $q = p$):

- If $t \le 1 - \pi_{\hat{A},k+\ell,\ell}^{(p)}k^{1-1/p}/\ell$, $\lambda \ge 4(2-t)t^{-1}(\sigma a\sqrt{2\ln(2d/\delta)/n})$, then $\|\Delta\hat{\beta}\|_p \le \frac{8}{t\omega_{\hat{A},k+\ell}^{(p)}}k^{1/p}\lambda$;

- If $t \le 1 - \gamma_{\hat{A},k+\ell,\ell}^{(p)}k^{1-1/p}/\ell$, $\lambda \ge 4(2-t)t^{-1}(\sigma a\sqrt{2\ln(2d/\delta)/n})$, then $\|\Delta\hat{\beta}\|_p \le \frac{8}{t}\lambda(k+\ell)^{1/p}/\mu_{\hat{A},k+\ell}^{(p)}$.

That is, if there exist $\ell \ge k$, $t \in (0,1)$, and $p \in [1,\infty]$ so that:

- $\lambda \ge 4(2-t)t^{-1}(\sigma a\sqrt{2\ln(2d/\delta)/n})$;
- either $t \le 1 - \pi_{\hat{A},k+\ell,\ell}^{(p)}k^{1-1/p}/\ell$, $\alpha \ge 8\varepsilon^{-1}(t\omega_{\hat{A},k+\ell}^{(p)})^{-1}k^{1/p}\lambda$; or $t \le 1 - \gamma_{\hat{A},k+\ell,\ell}^{(p)}k^{1-1/p}/\ell$, $\alpha \ge 8\varepsilon^{-1}(t\mu_{\hat{A},k+\ell}^{(p)})^{-1}(k+\ell)^{1/p}\lambda$;

then $\|\Delta\hat{\beta}\|_\infty \le \|\Delta\hat{\beta}\|_p \le \varepsilon\alpha$. Note that the condition $\|\Delta\hat{\beta}\|_\infty \le \varepsilon\alpha$ implies that $\mathrm{supp}_{(1+\varepsilon)\alpha}(\bar{\beta}) \subset \mathrm{supp}_\alpha(\hat{\beta}) \subset \mathrm{supp}_{(1-\varepsilon)\alpha}(\bar{\beta})$. This proves the desired result.

10.10. *Proof of Corollary 7.1.* We take $p = \infty$, $\ell = k$ and $t = 0.5$ in Theorem 7.1. Proposition 3.2 implies that $\omega_{\hat{A},k+\ell}^{(\infty)} \ge 1 - M_{\hat{A}}(k+\ell) \ge 0.5$ and $\pi_{\hat{A},k+\ell,\ell}^{(\infty)} \le M_{\hat{A}}\ell/(1 - M_{\hat{A}}(k+\ell)) \le 0.5$ under the conditions of the corollary. Therefore, the condition $8(\varepsilon\alpha\omega_{\hat{A},k+\ell}^{(\infty)})^{-1}\lambda \le t \le 1 - \pi_{\hat{A},k+\ell,\ell}^{(\infty)}k/\ell$ in Theorem 7.1 holds.

10.11. *Proof of Corollary 7.2.* We want to apply Theorem 7.1 with $\varepsilon = 0.5$, $\hat{A} = \hat{A}_n$, $\delta = \delta_n = \exp(-n^{s'})$, $t = 0.5$, $\alpha = \alpha_n = n^{-s/2}$, $d = d_n$, $k = k_n$, and $\ell = q_n k_n$. Then, $\lambda = \lambda_n = 4(2-t)t^{-1}(\sigma a\sqrt{2\ln(2d_n/\delta_n)/n})$. The condition $\rho_{\hat{A}_n,(1+2q_n)k_n}^{(2)} \le (1+q_n)\mu_{\hat{A}_n,(1+2q_n)k_n}^{(2)}$ implies that $\ell/k \ge \rho_{\hat{A},k+2\ell}^{(2)}/\mu_{\hat{A},k+2\ell}^{(2)} - 1 \ge 4(\pi_{\hat{A},k+\ell,\ell}^{(2)})^2\ell^{-1}$, where the second inequality is due to Proposition 3.1. Therefore, the condition $t \le 1 - \pi_{\hat{A},k+\ell,\ell}^{(2)}k^{0.5}/\ell$ is satisfied.

Now, from the conditions $k_n = o(n^{1-s-s'})$ and $k_n \ln(d_n) = o(n^{1-s})$, we have $\lim_{n\to\infty}(\sqrt{k_n}\lambda_n/\alpha_n) = 0$. Since $(\mu^{(2)}_{\hat{A},k+\ell})^{-1} = O(1)$, the condition $\alpha \geq 16(t\mu^{(2)}_{\hat{A},k+\ell})^{-1}\sqrt{k}\lambda$ is also satisfied when $n$ is sufficiently large.

Therefore, by Theorem 7.1, when $n$ is sufficiently large, $\mathrm{supp}_{1.5\alpha_n}(\bar{\beta}_n) \subset \mathrm{supp}_{\alpha_n}(\hat{\beta}_n) \subset \mathrm{supp}_{0.5\alpha_n}(\bar{\beta}_n)$ with probability at least $1 - \delta_n$. Since $1/\min_{j\in\mathrm{supp}_0(\bar{\beta}_n)}|\bar{\beta}_{n,j}| = o(n^{s/2})$, we know when $n$ is sufficiently large, $\min_{j\in\mathrm{supp}_0(\bar{\beta}_n)}|\bar{\beta}_{n,j}| > 2\alpha_n$. Thus, $\mathrm{supp}_{1.5\alpha_n}(\bar{\beta}_n) = \mathrm{supp}_{0.5\alpha_n}(\bar{\beta}_n) = \mathrm{supp}_0(\bar{\beta}_n)$. This means that $\mathrm{supp}_0(\bar{\beta}_n) = \mathrm{supp}_{\alpha_n}(\hat{\beta}_n)$ with probability at least $1 - \delta_n$ when $n$ is sufficiently large.

10.12. *Proof of Theorem 8.1.* Let $F = \mathrm{supp}_\lambda(\bar{\beta})$. We would like to apply Lemmas 10.2, 10.4 and 10.3, with $\hat{F} = \mathrm{supp}_\alpha(\hat{\beta})$ and $F_0 = \mathrm{supp}_{1.5\alpha}(\bar{\beta})$.

First, Theorem 7.1 implies that, with probability larger than $1 - \delta$, $\mathrm{supp}_{1.5\alpha}(\bar{\beta}) \subset \mathrm{supp}_\alpha(\hat{\beta}) \subset \mathrm{supp}_{0.5\alpha}(\bar{\beta}) \subset F$. Moreover,

$$\|\hat{\varepsilon}\|_\infty \leq a\sigma\sqrt{\frac{2\ln(2d/\delta)}{n}}. \tag{7}$$

This can be seen from the proof of Theorem 4.1 (which directly implies Theorem 7.1).

Let $\mathbf{z}_i = [\mathbf{z}_{i,1}, \ldots, \mathbf{z}_{i,d}] \in R^d$, so that $\mathbf{z}_{i,j} = -\frac{1}{n}\mathbf{x}_{i,j}$ if $j \in F_0$ and $\mathbf{z}_{i,j} = 0$, otherwise. We thus have $\hat{\varepsilon}_{F_0} = \sum_{i=1}^n (\mathbf{y}_i - \mathbf{E}\mathbf{y}_i)\mathbf{z}_i$. Since each $\mathbf{y}_i - \mathbf{E}\mathbf{y}_i$ is an independent sub-Gaussian random variable, and $\sum_{i=1}^n \|\mathbf{z}_i\|_2^2 = \sum_{j\in F_0}\sum_{i=1}^n (\mathbf{x}_{i,j}/n)^2 \leq qa^2/n$, we obtain, from Proposition 10.2, that, with probability larger than $1 - \delta$,

$$\|\hat{\varepsilon}_{F_0}\|_2 \leq a\sigma(1 + \sqrt{20\ln(1/\delta)})\sqrt{q/n}. \tag{8}$$

Therefore, with probability exceeding $1 - 2\delta$, both (7) and (8) hold. Therefore, Lemma 10.4 implies that

$$\begin{aligned}
\|\Delta\hat{\beta}'_I\|_2 &\leq \frac{2}{t\mu^{(2)}_{\hat{A},k+\ell}}[4\theta^{(2)}_{\hat{A},k+\ell,\ell}\ell^{-1}\|\bar{\beta}_G\|_1 + \rho^{(2)}_{\hat{A},k+\ell}\|\bar{\beta}_G\|_2 + \sqrt{k-q}\lambda + \|\hat{\varepsilon}_{F_0}\|_2] \\
&\quad + \|\bar{\beta}_G\|_2 \\
&\leq \frac{2}{t\mu^{(2)}_{\hat{A},k+\ell}}[5\rho^{(2)}_{\hat{A},k+\ell}\|\bar{\beta}_G\|_2 + \sqrt{k-q}\lambda + a\sigma(1 + \sqrt{20\ln(1/\delta)})\sqrt{q/n}] \\
&\quad + \|\bar{\beta}_G\|_2.
\end{aligned}$$

In the first inequality, we have used $\pi^{(2)}_{\hat{A},k+\ell,\ell} \leq \theta^{(2)}_{\hat{A},k+\ell,\ell}/\mu^{(2)}_{\hat{A},k+\ell}$ (Proposition 3.1). In the second inequality, we have used $\|\bar{\beta}_G\|_1 \leq \sqrt{s}\|\bar{\beta}_G\|_2 \leq \sqrt{\ell}\|\bar{\beta}_G\|_2$, and $\theta^{(2)}_{\hat{A},k+\ell,\ell} \leq \rho^{(2)}_{\hat{A},k+\ell}\sqrt{\ell}$ (Proposition 3.1). By combining the above estimate with Lemma 10.3, we obtain the desired bound.

10.13.  *Proof of Corollary 8.1.*  We take $p = \infty$, $\ell = s$, and $t = 0.5$ in Theorem 8.1. Proposition 3.2 implies that $\pi^{(2)}_{\hat{A},k+\ell,\ell} \leq M_{\hat{A}}(2s)^{1/2}s/(1-M_{\hat{A}}2s) \leq \sqrt{2s}/4$, $\omega^{(p)}_{\hat{A},s+\ell} \geq 1 - M_{\hat{A}}2s \geq 2/3$ and $\pi^{(p)}_{\hat{A},s+\ell,\ell} \leq M_{\hat{A}}s/(1-M_{\hat{A}}2s) \leq 1/4$.

Now, it is clear that $t \leq 1 - \pi^{(2)}_{\hat{A},k+\ell,\ell}k^{0.5}/\ell$ holds. Moreover, the condition $16(\alpha\omega^{(p)}_{\hat{A},s+\ell})^{-1}s^{1/p}\lambda \leq t \leq 1 - \pi^{(p)}_{\hat{A},s+\ell,\ell}s^{1-1/p}/\ell$ is also valid. Therefore, Theorem 8.1 can be applied with $\mu^{(2)}_{\hat{A},k+\ell} \geq 1 - M_{\hat{A}}(k+\ell) \geq 2/3$ and $\rho^{(2)}_{\hat{A},k+\ell} \leq 1 + M_{\hat{A}}(k+\ell) \leq 4/3$.

## 11. Conclusion.

This paper considers the performance of least squares regression with $L_1$ regularization from parameter estimation accuracy and feature selection quality perspectives. To this end, a general theorem is established in Section 4.

An important consequence of this theorem is a performance bound for Lasso similar to that of [4] for the Dantzig selector. The detailed comparison is given in Section 6. Our result gives an affirmative answer to an open question in [10] concerning whether a bound similar to that of [4] holds for Lasso. Another important consequence of Theorem 4.1 is the feature selection quality of Lasso using a nonzero thresholding feature selection method, which extends the zero thresholding method considered in [20]. Our method can remove some limitations of [20], as discussed in Section 7.

Moreover, we pointed out that the standard (one-stage) Lasso may be suboptimal under certain conditions. However, the problem can be remedied by combining the parameter estimation and feature selection perspectives of Lasso. In Section 8, a two-stage $L_1$-regularization procedure with selective penalization was analyzed. In practice, if one is able to appropriately tune the thresholding parameter using cross-validation, then the procedure should not be much worse than the standard one-stage Lasso. Theoretically, it is shown that, if the target vector can be decomposed as the sum of a sparse parameter vector with large coefficients and another (less sparse) vector with small coefficients, then the two-stage $L_1$-regularization procedure can lead to improved performance when $d$ is large.

Finally, we shall point out some limitations of our analysis. First, procedures considered in this work are not adaptive. For example, in the one-stage method, the regularization parameter $\lambda$ has to satisfy certain conditions that depend on $t$ and the noise level $\sigma$. In feature selection and the two-stage method, the threshold parameter $\alpha$ also needs to satisfy certain conditions.

Although, in practice, such parameters can be tuned using cross-validation, it still remains an interesting problem to come out with a theoretical procedure for setting such parameters that leads to a so-called "adaptive" estimation method. Moreover, although bounds in this paper can be applied in random design situations with small modifications, the results are incomplete for random design because the conditions on the design matrix (which is now random) needs to be shown to concentrate at a certain rate. Although a number of such results exist in the random matrix literature, a more general treatment with better integration is still needed in future work.

**Acknowledgments.** The author would like to thank Trevor Hastie, Terence Tao and an anonymous reviewer for helpful remarks that clarified differences between $L_1$ regularization and the Dantzig selector.

STATISTICS DEPARTMENT
RUTGERS UNIVERSITY
PISCATAWAY, NEW JERSEY 08854
USA
E-MAIL: tzhang@stat.rutgers.edu